How this translation was made:
The starting point was a recommendation by Urs Frauenfelder, and the collected works of Lyapunov at http://books.e-heritage.ru/book/10084104 which I captured in PDF. Based on the date (1896), it is clear that the document was created by traditional type-setting and later scanned to create the PDF. Then I went through the following steps:

- Opened the PDF document in Microsoft Word, which automatically ran OCR (optical character recognition), and converted it to an editable Word document.
- Corrected obvious errors, e.g. due to hyphenation and line breaks, and to formulas inline in the text.
- Recreated all formulas and equations using the built-in Microsoft Equation Editor.
- Translated the Russian text to English using the built-in Microsoft Translator.
- Edited the text to apply relevant mathematics conventions, e.g. referring to Hill's series instead of Hill's rows.

Finally, Marina Klein and Tanya Znamenskaya proofread the whole document to insure correct translation of Russian to English.

Acknowledgment
I would like to thank Urs Frauenfelder, Marina Klein, and Tanya Znamenskaya for their valuable support.

# ON THE SERIES PROPOSED BY HILL FOR THE REPRESENTATION OF THE MOTION OF THE MOON

1. In the first volume of the «American Journal of Mathematics» in 1878 G. Hill published a very interesting memoir «Researches in the Lunar Theory»[1] where the motion of the moon under the influence of the earth and the sun is interpreted as the limiting case of the motion of three bodies, to which it is possible to move, assuming the mass of the sun and its distance from the earth increase infinitely, and so that the ratio of the cube of this distance to the mass of the sun is closer to a certain limit independent of time. The solution to the problem in this setting, which, in the name of further study of the issue, Hill considers only as preliminary, delivers such an approximation, which takes into account the average movement of the sun, but neglects its parallax and eccentricity. For most of his findings, Hill also neglects the inclination of the lunar orbit and, assuming that the moon moves in the plane of the ecliptic, for its differential equations, by which it determines the relative movement of the moon orbiting the earth together with the sun, seeks not a common integral, but some special solution, subordinated to the condition of periodicity. This special solution contains two arbitrary constants, one of the connections and the average movement of the moon, and in Hill's further research this was supposed to play the role of a first approximation, and it is mainly studied in the memoir under consideration.[2]

Hill's special solution defines trigonometric series, in which the coefficients, known to depend on the average movements of the sun and the moon, are calculated by him by successive approximations. But when it comes to calculating these coefficients, Hill is not at all concerned with the convergence of his series.

Some general considerations suggest that the series in question should undoubtably converge if the ratio of the average movement of the sun to the average movement of the moon does not exceed a certain limit, and in this case, the series represent one of the periodic movements, the existence of which was discovered by Poincaré in his general research on periodic solutions to the differential equations of the three-body problem. The limit, as far as I know, was not sought by anyone, and it has not yet been proven that it is large enough for the convergence of the series under consideration at the magnitude of the average movements of the sun and moon, accepted in astronomy.

Meanwhile, given the importance that Hill's series may have in lunar theory, the question of convergence of them deserves serious attention. I found it therefore not unhelpful to publish my research on the issue, which I suggest in this article. The study shows that the magnitude of the relationship between the movements of the sun and the moon, taken in astronomy, quite ensures the convergence of the Hill series.

2. Hill's differential equations of moon motion are the following:

$$\left.\begin{array}{l}\dfrac{d^2x}{dt^2} - 2n\dfrac{dy}{dt} + \dfrac{\mu x}{(x^2+y^2)^{\frac{3}{2}}} = 3n^2x,\\[2mm] \dfrac{d^2y}{dt^2} + 2n\dfrac{dx}{dt} + \dfrac{\mu y}{(x^2+y^2)^{\frac{3}{2}}} = 0.\end{array}\right\} \quad (1)$$

Here $x$ and $y$ are the rectangular coordinates of the moon in relation to the axes, the beginning of which coincides with the center of the earth and which move in the plane of the ecliptic so that the $x$ axis is constantly directed from the earth to the sun, and the $y$ axis is in the direction of the year-long movement of the sun relative to the earth; $n$ means the average movement of the sun and $\mu$ the sum of the masses of the earth and the moon, multiplied by the ratio of attraction.

---

[1] A very detailed analysis of this memoir can be found in the third volume of the essay: Tisserand, Traité Mécanique Céleste.

[2] Recently, «Astronomical Journal» (vol. XV, No. 353) Hill's new memoir about the moon's movement, in which Hill examines a similar periodic decision, taking into account the parallax of the sun.

With these equations, Hill tries to satisfy the series equations

$$x = A_0 \cos n_1(t - t_0) + A_1 \cos 3n_1(t - t_0) + A_2 \cos 5n_1(t - t_0) + \cdots$$
$$y = B_0 \sin n_1(t - t_0) + B_1 \sin 3n_1(t - t_0) + B_2 \sin 5n_1(t - t_0) + \cdots$$

where $t_0$ and $n_1$ are arbitrary constants, and $A_0, A_1$ etc. $B_0, B_1$ etc. are constants, defined as functions of $n$ and $n_1$, from the condition of the task.

The constant $n_1$ represents the synodal average movement of the moon.

Assuming

$$n_1(t - t_0) = \tau$$

and, in general, for all $s$

$$A_s = a_s + a_{-s-1}, \quad B_s = a_s - a_{-s-1},$$

Hill presents his series in the form of

$$x = \sum a_s \cos(2s + 1)\tau, \quad y = \sum a_s \sin(2s + 1)\tau,$$

where the summation extends over all whole values $s$ from $-\infty$ to $+\infty$.

Instead of variables $x$ and $y$ Hill introduces then the following:

$$u = x + y\sqrt{-1}, \quad v = x - y\sqrt{-1},$$

for which the previous series gives

$$u = \sum a_s e^{(2s+1)i\tau}, \quad v = \sum a_{-s-1} e^{(2s+1)i\tau}, \qquad (2)$$

where $i = \sqrt{-1}$, and the differential equations, if we put

$$\frac{n}{n_1} = m, \quad \frac{\mu}{n_1^2} = k,$$

look like this:

$$\left.\begin{aligned}\frac{d^2u}{d\tau^2} + 2mi\frac{du}{d\tau} + \frac{ku}{(uv)^{\frac{3}{2}}} &= \frac{3}{2}m^2(u + v), \\ \frac{d^2v}{d\tau^2} - 2mi\frac{dv}{d\tau} + \frac{kv}{(uv)^{\frac{3}{2}}} &= \frac{3}{2}m^2(u + v).\end{aligned}\right\} \qquad (3)$$

Denoting by $C$ an arbitrary constant, from the equations we obtain:

$$\frac{du}{d\tau}\frac{dv}{d\tau} - \frac{2k}{\sqrt{uv}} = +\frac{3}{4}m^2(u + v)^2 - C, \qquad (4)$$

representing the force equation for the task at hand.

Excluding from the equations (3) and (4) the constant $k$, we obtain

$$\left.\begin{aligned}v\frac{d^2u}{d\tau^2} - u\frac{d^2v}{d\tau^2} + 2mi\frac{d(uv)}{d\tau} + \frac{3}{2}m^2(u^2 - v^2) &= 0, \\ v\frac{d^2u}{d\tau^2} + u\frac{d^2v}{d\tau^2} + \frac{du}{d\tau}\frac{dv}{d\tau} + 2mi\left(v\frac{du}{d\tau} - u\frac{dv}{d\tau}\right) - \frac{9}{4}m^2(u + v)^2 + C &= 0,\end{aligned}\right\} \qquad (5)$$

To satisfy the equations (3), Hill first tries to satisfy the equations (5), and arbitrary functions $u$ and $v$, satisfying these latest equations, the following need to be satisfied:

$$\frac{d^2u}{d\tau^2} + 2mi\frac{du}{d\tau} + \frac{hu}{(uv)^{\frac{3}{2}}} = \frac{3}{2}m^2(u+v),$$

$$\frac{d^2v}{d\tau^2} - 2mi\frac{dv}{d\tau} + \frac{hv}{(uv)^{\frac{3}{2}}} = \frac{3}{2}m^2(u+v),$$

where $h$ is some constant, which for the final solution of the problem it remains only to satisfy any ratio between the constant, leading to equality

$$h = k. \tag{6}$$

Looking for the conditions under which the series (2) satisfy the equations (5), we get the expression of constant $C$ for the coefficients $a_s$, and the following ratios between the latter:

$$4j\sum_s (2s - j + 1 + m)a_s a_{s-j} - \frac{3}{2}m^2 \sum_s a_s(a_{j-s-1} - a_{-j-s-1}) = 0,$$

$$\sum_s \left[(2s+1)(2s-2j+1) + 4j^2 + 4(2s-j+1)m + \frac{9}{2}m^2\right]a_s a_{s-j}$$
$$+ \frac{9}{4}m^2 \sum_s a_s(a_{j-s-1} - a_{-j-s-1}) = 0,$$

which lead to

$$\sum_s [8s^2 - 8(4j-1)s + 20j^2 - 16j + 2 + 4(4s - 5j + 2)m + 9m^2]a_s a_{s-j} + 9m^2 \sum_s a_s a_{j-s-1} = 0,$$

$$\sum_s [8s^2 + 8(2j+1)s - 4j^2 + 8j + 2 + 4(4s + j + 2)m + 9m^2]a_s a_{s-j} + 9m^2 \sum_s a_s a_{-j-s-1} = 0,$$

Here, the summation extends to all whole values $s$ from $-\infty$ to $+\infty$, and $j$ can have all whole values except zero. But by giving both positive and negative values, one can limit oneself to considering only one of the two equations written now, for it is not difficult to see that one is derived from the other by substituting $j$ for $-j$.

Hill examines a combination of these equations.

If the first of them is presented by equality $L_j = 0$ and, therefore, the second equality is $L_{-j} = 0$, that combination this can be formulated as:

$$[4j^2 - 8j - 2 - 4(j+2)m - 9m^2]L_j + [20j^2 - 16j + 2 - 4(5j - 2)m + 9m^2]L_{-j} = 0.$$

The resulting equation is easily brought to the mind

$$\sum_s ([j,s]a_s a_{s-j} + [j]a_s a_{j-s-1} + (j)a_s a_{-j-s-1}) = 0, \tag{7}$$

where, according to Hill's designation,

$$[j,s] = -\frac{s}{j}\frac{4s(j-1) + 4j^2 + 4j - 2 - 4(s - j + 1)m + m^2}{2(4j^2 - 1) - 4m + m^2},$$

$$[j] = -\frac{3m^2}{16j^2}\frac{4j^2 - 8j - 2 - 4(j+2)m - 9m^2}{2(4j^2 - 1) - 4m + m^2},$$

$$(j) = -\frac{3m^2}{16j^2}\frac{20j^2 - 16j + 2 - 4(5j - 2)m + 9m^2}{2(4j^2 - 1) - 4m + m^2}.$$

In equation (7) Hill makes use of the relationships $\frac{a_s}{a_0}$ ($s = \pm 1, \pm 2, \pm 3, ...$), admitting that these relationships are functions of the parameter $m$ and disappear with $m$ so that for $m$ infinitesimally small first order $\frac{a_s}{a_0}$ there's an infinitesimal order $\pm 2s$. The allows us to determine the relationship in question by successive approximations with precision to members of any order of $m$, which finally leads to the representation of these relationships in the form of endless series.

Let's point out here the series to which you can come from considering the method of successive approximations used by Hill.

Hill's calculations are based on the fact that in the equation (7) the coefficients $[j]$ and $(j)$ disappear for $m = 0$ and the second order of $m$, and the coefficients $[j, s]$ such that

$$[j, j] = -1, \quad [j, 0] = 0. \tag{8}$$

Highlighting the factor $m^2$, we put

$$[j] = \overline{[j]}m^2, \quad (j) = \overline{(j)}m^2$$

and then by introducing the parameter $\lambda$, equation (7) is replaced by:

$$\sum_s ([j,s]a_s a_{s-j} + \overline{[j]}\lambda a_s a_{j-s-1} + \overline{(j)}\lambda a_s a_{-j-s-1}) = 0. \tag{9}$$

Now assuming

$$\frac{a_s}{a_0} = \lambda^{|s|}\alpha_s, \tag{10}$$

where $|s|$ denotes the numerical value $s$, let's look for $\alpha_s$ ($s = \pm 1, \pm 2, \pm 3, ...$) in the form of series, arranged by whole positive powers of $\lambda$.

Substituting these series into expression (10) and then replacing $\lambda$ by $m^2$, and getting these series, which leads to the way Hill calculated them, if only members of the orders of which are below the margin of error are retained.

It is easy to see that the series that represent the magnitudes $\alpha_s$, will only contain even powers of $\lambda$ and that the coefficients in these series will be determined in a known sequence so that the final number of algebraic actions will suffice to determine each factor in each series.

In fact, after the situation (10) in the equation (9) will meet only those powers of $\lambda$, indicators of which are in the form of

$$|s| + |s - j|, \quad |s| + |s - j + 1| + 1, \quad |s| + |s + j + 1| + 1,$$

and the latter obviously give only numbers equal $|j|$ or superior to $|j|$ on the even number.

So, if after a substitution (10) all members of the equation (9) divide into $\lambda^{|j|}$ it will meet only positive and even powers of $\lambda$, and because of (8) it will be possible to give the following form:

$$\alpha_j = L_{j,0} + \sum_{s=1}^{\infty} L_{j,s}\lambda^{2s},$$

where

$$L_{1,0} = \overline{[1]}, \quad L_{-1,0} = \overline{(-1)}$$

and in general

$$L_{j,0} = \begin{cases} \overline{[j]}\sum_{s=0}^{j-1} \alpha_s \alpha_{j-s-1} + [j,s]\sum_{s=1}^{j-1} \alpha_s \alpha_{s-j}, & \text{if } j > 1, \\ \overline{(j)}\sum_{s=0}^{-j-1} \alpha_s \alpha_{-j-s-1} + [j,s]\sum_{s=j+1}^{-1} \alpha_s \alpha_{s-j}, & \text{if } j < -1, \end{cases}$$

and all $L_{j,s}$ there is a set of members of the species

$$C\alpha_p\alpha_q,$$

for which symbols $p$ and $q$ satisfy the equality

$$p - q = j, \quad |p| + |q| = 2s + |j|$$

or

$$p + q = \pm j - 1, \quad |p| + |q| = 2s + |j| - 1.$$

It is clear from this that the series under consideration will be a kind of

$$\alpha_s = \alpha_{s,0} + \alpha_{s,1}\lambda^2 + \alpha_{s,2}\lambda^4 + \cdots$$

and that coefficients $\alpha_{s,j}$ can be put in a linear series, in which each member will be determined by the previous one, and the first member will be some of the coefficients $\alpha_{1,0}, \alpha_{-1,0}$ defined directly by the formulas

$$\alpha_{1,0} = \overline{[1]}, \quad \alpha_{-1,0} = \overline{(-1)}.$$

You can see that each of the coefficients $\alpha_{s,j}$, will be a rational function of $m$, in which, if it is presented under the form of an unreduced fraction, the denominator will contain only the multipliers of the type

$$2(4\sigma^2 - 1) - 4m + m^2, (\sigma = 1,2,3,\dots),$$

which is not equal to zero for any real physical values of $m$

Defining the ratio $\frac{a_s}{a_0}$, Hill moves on to finding $a_0$ and composing in a variety of ways equation (6), gets a few different expressions for $a_0$. One of these expressions, deduced from the consideration of members of the first power of $e^{i\tau}$ as a result of the substitution in series (2) in the first of the equations (3), we have the following:

$$a_0 = \frac{(kJ)^{1/3}}{\left(1 + 2m + \frac{3}{2}m^2 + \frac{3}{2}m^2 \frac{a_{-1}}{a_0}\right)^{1/3}}, \tag{11}$$

where $J$ means not dependent on $\zeta$ member in decomposition expression

$$\left\{\sum \frac{a_s}{a_0} \zeta^{2s}\right\}^{-1/2} \left\{\sum \frac{a_s}{a_0} \zeta^{-2s}\right\}^{-3/2}$$

by positive and negative powers of $\zeta$.

Applying these formulas to lunar theory, Hill takes

$$m = 0.08084\ 89338\ 08312$$

and at the final conclusions neglects only the members above 13-th power relative to $m$. Computing at the same time the ratio $\frac{a_s}{a_0}$ with fifteen decimal places, he assumes that the error does not exceed the last decimal point. This is the conclusion, however, based only on a comparison of calculated members of different orders and cannot be considered proven, because in fact Hill does not give any means to determine the higher margin of error in calculating the ratios $\frac{a_s}{a_0}$ in the way he proposes.

Setting out the essence of Hill's analysis, let's move on to the question of the convergence of the series, that will be determined, however, a bit differently.

3. It's easy to see that equation (9), which, entering parameter $\lambda$, we replaced the above equation (7), it's possible to find a series of the form (2), satisfying the following differential equations:

$$\left.\begin{array}{l} \dfrac{d^2u}{d\tau^2} + 2mi\dfrac{du}{d\tau} - \dfrac{3}{2}m^2 u + \dfrac{ku}{(uv)^{\frac{3}{2}}} = \dfrac{3}{2}\lambda v, \\[2mm] \dfrac{d^2v}{d\tau^2} - 2mi\dfrac{dv}{d\tau} - \dfrac{3}{2}m^2 v + \dfrac{kv}{(uv)^{\frac{3}{2}}} = \dfrac{3}{2}\lambda u, \end{array}\right\} \tag{12}$$

that, by conversion to the variables $x$ and $y$, takes the form

$$\left.\begin{array}{l}\dfrac{d^2x}{d\tau^2} - 2m\dfrac{dy}{d\tau} - \dfrac{3}{2}m^2x + \dfrac{kx}{(x^2+y^2)^{\frac{3}{2}}} = \dfrac{3}{2}\lambda x, \\ \dfrac{d^2y}{d\tau^2} + 2m\dfrac{dx}{d\tau} - \dfrac{3}{2}m^2y + \dfrac{ky}{(x^2+y^2)^{\frac{3}{2}}} = -\dfrac{3}{2}\lambda y. \end{array}\right\} \quad (13)$$

We're now going to look at these equations instead of equations (1), and let's keep in mind that in the final result we have $\lambda = m^2$.

Putting into equations (13) $\lambda = 0$, we get equations, which can be satisfied, by taking

$$x = a\cos\tau, \qquad y = a\sin\tau \quad (14)$$

and using for $a$ a constant defined by

$$\frac{k}{a^3} = 1 + 2m + \frac{3}{2}m^2.$$

This constant will represent the length, a little different from the average distance of the moon from the earth, the cube of which is equal to

$$\frac{k}{(1+m)^2}.$$

Going back now to the equations (13), in which $\lambda$ will be assumed numerically not exceeding a certain limit, we'll look for a solution for them, in which $x$ and $y$ will be periodic functions of $\tau$ with period $2\pi$, when applying $\lambda = 0$ in expression (14).

We will thus seek periodic movement, which at small values of $\lambda$ would be little different from uniform movement in a circle of radius $a$ centered in the center of the earth.

We can present the solution with formulas

$$x = a(1+\xi)\cos\tau - a\eta\sin\tau, \quad y = a(1+\xi)\sin\tau + a\eta\cos\tau,$$

considering $\xi$ and $\eta$, periodic function $\tau$, and disappear at $\lambda = 0$.

Instead of these features, we will further consider the following combinations:

$$-\xi - i\eta = p, \quad -\xi + i\eta = q,$$

by means of the variables introduced above $u$ and $v$ are expressed by formulas

$$u = a(1-p)e^{i\tau}, \quad v = a(1-q)e^{-i\tau},$$

and by virtue of the latest equations (12) will deliver the following differential equations to determine the functions $p$ and $q$:

$$\frac{d^2p}{d\tau^2} + 2(1+m)i\frac{dp}{d\tau} + l(1-p) - \frac{l}{(1-p)^{1/2}(1-q)^{3/2}} = \frac{3}{2}\lambda e^{-2i\tau},$$

$$\frac{d^2q}{d\tau^2} - 2(1+m)i\frac{dq}{d\tau} + l(1-q) - \frac{l}{(1-p)^{3/2}(1-q)^{1/2}} = \frac{3}{2}\lambda e^{+2i\tau},$$

where we reduce the length by using

$$l = 1 + 2m + \frac{3}{2}m^2.$$

The features included here

$$(1-p)^{-1/2}(1-q)^{-3/2}, \quad (1-p)^{-3/2}(1-q)^{-1/2}$$

will represent the series located on whole positive values of $p$ and $q$, so our analysis will only apply to the case when the module functions $p$ and $q$ always remain smaller than one.

Meaning that in the named series the aggregate elements above the first power relative to $p$ and $q$ via $R_p$ and $R_q$, we have

$$(1-p)^{-1/2}(1-q)^{-3/2} = 1 + \frac{1}{2}p + \frac{3}{2}q + R_p$$

$$(1-p)^{-3/2}(1-q)^{-1/2} = 1 + \frac{3}{2}p + \frac{1}{2}q + R_q$$

and our equations will be:

$$\left. \begin{aligned} \frac{d^2p}{d\tau^2} + 2(1+m)i\frac{dp}{d\tau} - \frac{3}{2}l(1+p) &= \frac{3}{2}\lambda(q-1)e^{-2i\tau} + lR_p, \\ \frac{d^2q}{d\tau^2} - 2(1+m)i\frac{dq}{d\tau} - \frac{3}{2}l(1+p) &= \frac{3}{2}\lambda(p-1)e^{+2i\tau} + lR_q. \end{aligned} \right\} \quad (15)$$

For these equations, we will now look for periodic solutions with a period $2\pi$[3], with the requirement that the functions $p$ and $q$ for infinitely small $\lambda$ were also made infinitely small, and above all, let's show that whenever the module parameter $\lambda$ is small enough, the functions $p$ and $q$ in this solution can be determined by series, based on whole positive powers of $\lambda$, and that these series can only contain one arbitrary constant.

Let's use the Poincaré method for this.

Let $p_0, q_0, p'_0, q'_0$ be the initial value functions

$$p, \quad q, \quad p' = \frac{dp}{d\tau}, \quad q' = \frac{dq}{d\tau},$$

that, for example, correspond to $\tau = 0$.

Whenever modules of the magnitude

$$p_0, \quad q_0, \quad p'_0, \quad q'_0, \quad \lambda \quad (16)$$

are small enough, the functions $p$ and $q$, satisfy the equations (15), and their derivatives $p'$ and $q'$, based on one general statement, it is possible to represent series based on whole positive powers of (16), for all values of $\tau$, lying between 0 and some limit $T$, which is a decrease in module values (16) can be made as big as desired.

Suppose that in this fashion we have $T > 2\pi$. Then the now named series can be used to compile the values of functions $p, q, p', q'$, relevant to $\tau = 2\pi$.

Let

$$\bar{p}, \quad \bar{q}, \quad \bar{p}', \quad \bar{q}'$$

be these values.

Composing their equations

$$\bar{p} = p_0, \quad \bar{q} = q_0, \quad \bar{p}' = p'_0, \quad \bar{q}' = q'_0, \quad (17)$$

we'll get the conditions necessary and obviously sufficient for the functions $p$ and $q$, determined by equations (15), were periodic with a period of $2\pi$.

The question is thus brought to the study of equations (17), which must satisfy the initial condition constants $p_0, q_0, p'_0, q'_0$.

It's easy to show that equations (17) are not independent and that one of them is the consequence of the other three.

To do this, we notice that the equations (15) allow the following integral equation:

$$[p' - i(1-p)][q' + i(1-q)] - \frac{3}{2}m^2(1-p)(1-q) - \frac{2l}{\sqrt{(1-p)(1-q)}}$$
$$= \frac{3}{4}\lambda[(1-p)^2 e^{2i\tau} + (1-q)^2 e^{-2i\tau}] - \frac{C}{a^3},$$

---

[3] As we shall see, it's not just to address this $2\pi$ but also $\pi$ will be a period.

is the transformation of an equation that corresponds to (4) for equations (12).

Decomposing this way the function

$$(1-p)^{-1/2}(1-q)^{-1/2}$$

in the series by powers of $p$ and $q$, we bring this integral equation in the form

$$2(1+m)(p+q) - i(p'-q') + \cdots + \frac{3}{4}\lambda\left[(1-p)^2 e^{2i\tau} + (1-q)^2 e^{-2i\tau}\right] = \text{const},$$

where omitted members are above the second power of $p, q, p', q'$, and from here, substituting instead $p, q, p', q'$ from their expression as series, then putting $\tau = 2\pi$, we get the following identity:

$$2(1+m)(\bar{p}+\bar{q}) - i(\bar{p}'-\bar{q}') + \cdots + \frac{3}{4}\lambda\left[(1-\bar{p})^2 + (1-\bar{q})^2\right]$$
$$= 2(1+m)(p_0+q_0) - i(p_0'-q_0') + \cdots + \frac{3}{4}\lambda\left[(1-p_0)^2 + (1-q_0)^2\right]$$

relative to (16).

From here it is clear that if any three of the equations (17) are satisfied, then the fourth will be satisfied as well.

On this basis, we may limit ourselves to considering the following three equations:

$$\bar{p} = p_0, \quad \bar{p}' = p_0', \quad \bar{q}' = q_0'. \tag{18}$$

Let's see what these equations can give

Turning to equations (15), we notice that if you drop all the members above the first power relative to $p$ and $q$, as well as the members that depend on $\lambda$, those equations will turn into linear equations with constant coefficients, for which the common integral will appear in the form of

$$p = \frac{3}{2}lC_1 e^{-\chi i\tau} - \left[\chi^2 + 2(1+m)\chi + \frac{3}{2}l\right]C_2 e^{-\chi i\tau} + C_3\tau + C'$$

$$q = \frac{3}{2}lC_2 e^{-\chi i\tau} - \left[\chi^2 + 2(1+m)\chi + \frac{3}{2}l\right]C_1 e^{+\chi i\tau} - C_3\tau + C'',$$

where

$$\chi = \sqrt{1 + 2m - \frac{1}{2}m^2},$$

and $C_1, C_2, C_3, C', C''$ are arbitrary constants related by

$$C' + C'' = \frac{4}{3}\frac{m+1}{l}iC_3.$$

Expressing these constants in terms of the initial values of the functions $p, q, p', q'$, we'll present them in the form of some linear combinations of the values $p_0, q_0, p_0', q_0'$, and $C_1$ $C_2$ and $C_3$ will be expressed only with the help of the following three constants:

$$p_0 + q_0, \quad p_0', \quad q_0', \tag{19}$$

which will be independent linear forms.

Now let's consider

$$C_1, C_2, C_3 \quad \text{and} \quad C' - C'' = C$$

defined as linear forms of the values

$$p_0 + q_0, \quad p_0', \quad q_0' \quad \text{and} \quad p_0 - q_0,$$

we will introduce these forms as arbitrary constants and in the general integral of exact equations (15), which will be presented as series, based on whole positive powers of

$$C_1, C_2, C_3, C \quad \text{and} \quad \lambda.$$

We'll get with this for $p$ and $q$ expressions that will be different from those just considered by the members above the first power regarding $C_1, C_2, C_3, C$ and $\lambda$ members that depend on $\lambda$.

Because of that, making equations (18) and writing them only of members that are independent of $\lambda$ and linear relative to $C_1, C_2, C_3, C$, we find that these equations are of the form

$$\frac{3}{2}l(\rho-1)C_1 - \left[\chi^2 + 2(1+m)\chi + \frac{3}{2}l\right]\left(\frac{1}{\rho}-1\right)C_2 + 2\pi C_3 + \cdots = 0,$$

$$\frac{3}{2}l(\rho-1)C_1 + \left[\chi^2 + 2(1+m)\chi + \frac{3}{2}l\right]\left(\frac{1}{\rho}-1\right)C_2 + \cdots = 0,$$

$$\left[\chi^2 + 2(1+m)\chi + \frac{3}{2}l\right](\rho-1)C_1 + \frac{3}{2}l\left(\frac{1}{\rho}-1\right)C_2 + \cdots = 0,$$

where

$$\rho = e^{2\pi\chi i}$$

and because the functional dependence of the first parts of them in terms of the values $C_1, C_2, C_3$ after the position

$$C_1 = C_2 = C_3 = C = \lambda = 0$$

turns into the value

$$2\pi(\rho-1)\left(\frac{1}{\rho}-1\right)\left\{\frac{3}{4}l^2 - \left[\chi^2 + 2(1+m)\chi + \frac{3}{2}l\right]^2\right\},$$

or, equivalently,

$$-16\pi(1+m)\chi[\chi + 2(1+m)]^2\sin^2\pi\chi,$$

which, at the value of $m$, appropriate to lunar theory, is not zero, that based on the well-known statement of the equation, ours will be solvable in terms of $C_1, C_2, C_3$ and provided that these latter must disappear at $C = \lambda = 0$, give them quite certain expressions in the form of series, based on whole positive powers of $C$ and $\lambda$, absolutely converging, as far as the modules $C$ and $\lambda$ are small enough.

Because in the expressions that are obtained in this way $C_1, C_2$, and $C_3$ and the constant $C$ will only appear in the members above the first power relative to $C$ and $\lambda$, these expressions will allow us to determine the magnitude (19) in the form of series, based on whole positive powers of $p_0 - q_0$ and $\lambda$.

Thus, the initial values $p_0, q_0, p_0', q_0'$, and therefore also the functions $p, q, p', q'$, corresponding to the periodic solution, appear in series, based on whole positive powers of $p_0 - q_0$ and $\lambda$ and will contain one arbitrary constant $p_0 - q_0$.

We will further consider this periodic solution as an assumption

$$p_0 = q_0, \tag{20}$$

tantamount to the assumption that when $\tau = 0$ $y = 0$.

Because $\tau$ already contains an arbitrary constant $t_0$, that assumption would not impose any significant limitation on our problem and would only serve to determine what we would understand under $t_0$, which will present thus one of the moments of connections.

4. Based on the above, we will be looking for our periodic solution in the form of series

$$\left.\begin{array}{l} p = p_1\lambda + p_2\lambda^2 + p_3\lambda^3 + \cdots, \\ q = q_1\lambda + q_2\lambda^2 + q_3\lambda^3 + \cdots, \end{array}\right\} \tag{21}$$

based on whole positive powers of $\lambda$.

Here $p_j, q_j$, are independent of $\lambda$ periodic functions of $\tau$, which, according to the condition (20), will be determined by the assumption that

$$\text{for} \quad \tau = 0 \quad p_j = q_j. \tag{22}$$

Substituting the series (21) into the equation (15) and expressing that the latter are satisfied regardless of $\lambda$, first, we get the next system of equations:

$$\left.\begin{array}{l}\dfrac{d^2 p_1}{d\tau^2} + 2(1+m)i\dfrac{dp_1}{d\tau} - \dfrac{3}{2}l(p_1 + q_1) = -\dfrac{3}{2}e^{-2i\tau},\\ \dfrac{d^2 q_1}{d\tau^2} - 2(1+m)i\dfrac{dq_1}{d\tau} - \dfrac{3}{2}l(p_1 + q_1) = -\dfrac{3}{2}e^{+2i\tau},\end{array}\right\} \quad (23)$$

which will determine $p_1$ and $q_1$. Then, if we assume that decomposition by ascending powers of $\lambda$ after substitution of expressions (21) this gives

$$\left.\begin{array}{l}R_p = P_2 \lambda^2 + P_3 \lambda^3 + P_4 \lambda^4 \ldots,\\ R_q = Q_2 \lambda^2 + Q_3 \lambda^3 + Q_4 \lambda^4 \ldots,\end{array}\right\}$$

we get multiple systems of equations of the form

$$\left.\begin{array}{l}\dfrac{d^2 p_j}{d\tau^2} + 2(1+m)i\dfrac{dp_j}{d\tau} - \dfrac{3}{2}l(p_j + q_j) = \dfrac{3}{2}q_{j-1}e^{-2i\tau} + lP_j,\\ \dfrac{d^2 q_j}{d\tau^2} - 2(1+m)i\dfrac{dq_j}{d\tau} - \dfrac{3}{2}l(p_j + q_j) = \dfrac{3}{2}p_{j-1}e^{+2i\tau} + lQ_j,\end{array}\right\} \quad (24)$$

of which for $j = 2,3\ldots$ we can sequentially find $p_2$ and $q_2$, $p_3$ and $q_3$ etc., because $P_j$ and $Q_j$ will obviously depend only on the values

$$p_1, p_2, \ldots, p_{j-1}, q_1, q_2, \ldots, q_{j-1}. \quad (25)$$

These equations, subject to periodicity and the condition (22) functions $p_j, q_j$ will be determined and, as it is easy to see, there will be species

$$\left.\begin{array}{l}p_j = \displaystyle\sum_{s=0}^{s=j} a_{j,j-2s} e^{2(j-2s)i\tau},\\ q_j = \displaystyle\sum_{s=0}^{s=j} a_{j,j-2s} e^{-2(j-2s)i\tau},\end{array}\right\} \quad (26)$$

where $a_{j,\sigma}$, are essentially some constants.

Indeed, from the equations (23) under the said conditions, we find

$$p_1 = a_{1,1} e^{2i\tau} + a_{1,-1} e^{-2i\tau}, \quad q_1 = a_{1,1} e^{-2i\tau} + a_{1,-1} e^{2i\tau},$$

where

$$a_{1,1} = -\dfrac{9}{16}\dfrac{2 + 4m + 3m^2}{6 - 4m + m^2}, \quad a_{1,-1} = \dfrac{3}{16}\dfrac{38 + 28m + 9m^2}{6 - 4m + m^2},$$

and admitting that expressions (26) hold for all values of $j$, smaller values of $r$, it's easy to prove that they will hold for $j = r$.

To do this, we notice that $P_j$, and $Q_j$, as follows from their very definition, are polynomials made up of values (25) so that the weight of each member is equal to $j$. Therefore, if, for example, in terms of the expression

$$\dfrac{3}{2} q_{r-1} e^{-2i\tau} + lP_r$$

all belonging to $p_j$ and $q_j$ replace with the expressions (26), can be put into this form:

$$p_j = e^{2ji\tau} \sum_{s=0}^{s=j} a_{j,j-2s} e^{-4si\tau}, \quad q_j = e^{2ji\tau} \sum_{s=0}^{s=j} a_{j,2s-j} e^{-4si\tau},$$

we get this result

$$\frac{3}{2}q_{r-1}e^{-2i\tau} + lP_r = \sum_{s=0}^{s=j} A_{r,r-2s}e^{2(r-2s)i\tau},$$

where $A_{r,r-2s}$ is a constant dependent in a well-known way on the constant $a_{j,\sigma}$, relevant when $j < r$.

At the same time, we have

$$\frac{3}{2}p_{r-1}e^{2i\tau} + lQ_r = \sum_{s=0}^{s=j} A_{r,r-2s}e^{-2(r-2s)i\tau},$$

for $Q_r$ coming from $P_r$, by permutation of the letters $p$ and $q$, and in the formulas (26), $q_j$ is taken from $p_j$ by replacing $\tau$ by $-\tau$.

Consequently, if all $p_j$, and $q_j$, for which $j < r$, appear in expressions (26), the equation (24) for $j = r$ takes the form

$$\left.\begin{aligned}\frac{d^2p_r}{d\tau^2} + 2(1+m)i\frac{dp_r}{d\tau} - \frac{3}{2}l(p_r + q_r) &= \sum_{s=0}^{s=j} A_{r,r-2s}e^{2(r-2s)i\tau}, \\ \frac{d^2q_r}{d\tau^2} - 2(1+m)i\frac{dq_r}{d\tau} - \frac{3}{2}l(p_r + q_r) &= \sum_{s=0}^{s=j} A_{r,r-2s}e^{-2(r-2s)i\tau},\end{aligned}\right\}$$

and with the requirements we set, they will give $p_r$ and $q_r$ expressions of the same type as (26), and the coefficients $a_{r,\sigma}$ determined by equations of the form

$$\left[4\sigma^2 + 4(1+m)\sigma + \frac{3}{2}l\right]a_{r,\sigma} + \frac{3}{2}la_{r,-\sigma} = -A_{r,\sigma} \quad (\sigma = -r, -r+2, \ldots, r-2, r),$$

from which we find

$$a_{r,\sigma} = \frac{\frac{3}{2}lA_{r,-\sigma} - \left[4\sigma^2 + 4(1+m)\sigma + \frac{3}{2}l\right]A_{r,\sigma}}{2\sigma^2[2(4\sigma^2 - 1) - 4m + m^2]}, \tag{27}$$

if $\sigma$ is not zero and

$$a_{r,0} = -\frac{1}{3l}A_{r,0}. \tag{28}$$

The latter formula refers to the case of even $r$.

So we know that $p_1$ and $q_1$ are defined by the expressions of the type (26), we can be sure the same kind of expression will hold for all $p_j$, and $q_j$. At the same time, knowing $a_{1,1}$ and $a_{1,-1}$ we can calculate all the other coefficients by formulas (27) and (28), which give all $a_{r,\sigma}$, valid for all $r$, when all are known $a_{j,s}$, for which $j < r$.

In this way, we make sure that the functions we seek $p$ and $q$ will be presented in series

$$\left.\begin{aligned}p &= \sum_{j=1}^{j=\infty}\sum_{s=0}^{s=j} a_{j,j-2s}\lambda^j e^{2(j-2s)i\tau}, \\ q &= \sum_{j=1}^{j=\infty}\sum_{s=0}^{s=j} a_{j,j-2s}\lambda^j e^{2(2s-j)i\tau},\end{aligned}\right\} \tag{29}$$

Where to functions $\xi$ and $\eta$ we get the following expressions:

$$\left.\begin{array}{l}\xi = -\sum_{j=1}^{j=\infty}\sum_{s=0}^{s=j} a_{j,j-2s}\lambda^j \cos 2(j-2s)\tau, \\ \eta = -\sum_{j=1}^{j=\infty}\sum_{s=0}^{s=j} a_{j,j-2s}\lambda^j \sin 2(j-2s)\tau,\end{array}\right\} \quad (30)$$

Writing the series (29) for powers of $e^{2i\tau}$ and doing

$$a(1-p) = \sum a_s e^{2si\tau}, \quad a(1-q) = \sum a_s e^{-2si\tau},$$

we have

$$\left.\begin{array}{l}a_0 = a\left\{1 - \sum_{s=1}^{s=\infty} a_{2s,0}\lambda^{2s}\right\}, \\ a_{\pm\sigma} = -a\sum_{s=0}^{s=\infty} a_{\sigma+2s,\pm\sigma}\lambda^{\sigma+2s},\end{array}\right\} \quad (31)$$

where for $\sigma$ of course the number is positive.

Here for functions $u$ and $v$ we get the expressions (2), which, according to Hill's assumption, the coefficient $a_\sigma$ are such that when $\lambda$ the infinitely small first-order ratio $\frac{a_{\pm\sigma}}{a_0}$ there's an infinitesimally small order $\sigma$.

Formulas (31) give expressions of coefficients $a_s$ considered by Hill in the form of series, formulated by powers of $\lambda$. To get the same expressions from Hill's formulas, we would have to use formula (11) to replace the following:

$$a_0 = aJ^{1/3}\left\{1 + \frac{3}{2}\frac{\lambda}{l}\frac{a_{-1}}{a_0}\right\}^{-1/3} \quad (32)$$

and that expression for $a_0$, and then expressions

$$a_{\pm\sigma} = a_0\lambda^\sigma \sum_{s=0}^{s=\infty} \alpha_{\pm\sigma,s}\lambda^{2s} \quad (33)$$

to put in series by powers of $\lambda$.

This way we would get formulas, giving expressions of the coefficients $a_{r,\sigma}$ via coefficients $\alpha_{j,s}$, calculated by Hill formulas, as shown in part 2, and I have to see that, to calculate the formula, these are more convenient than ours. But it is very difficult to draw any conclusions about the series in question here, whereas our formulas (27) and (28) are very convenient for this purpose.

Turning now to the question of convergence, we will consider the series (29) as formulated by powers of $\lambda$ and $e^{2i\tau}$ and at the same time, limited to physical values of $\tau$, let's show that for the magnitude of $m$ and $\lambda$, corresponding to lunar theory, these series are absolutely convergent, from which will follow the convergence under the same conditions of the series (30) and (31).

We will see that, under said conditions, the series, composed of series member modules (29), will be less than 1. That's why inequality will still hold

$$\sum_{s=1}^{s=\infty} |a_{2s,0}|\lambda^{2s} < 1,$$

if, under $|a|$ we understand the absolute value of $a$.

Hence, due to the absolute convergence of the series (31), we conclude about the same convergence of series

$$\alpha_{\pm\sigma,0}\lambda^\sigma + \alpha_{\pm\sigma,1}\lambda^{\sigma+2} + \alpha_{\pm\sigma,2}\lambda^{\sigma+4} + \cdots,$$

that we have used the relationship $\frac{a_{\pm\sigma}}{a_0}$.

We will also look at the series that the series (29) turn into after the coefficients are replaced by $a_{j,r}$ their decomposition by powers of $m$.

These coefficients, as seen by the formulas (27) and (28), will be rational functions of $m$ and will consist of fractions, the denominators of which will only contain factors of the type

$$2(4\sigma^2 - 1) - 4m + m^2, \qquad (\sigma = 1,2,3,...)$$

and

$$l = 1 + 2m + \frac{3}{2}m^2.$$

That's why when $|m| < \sqrt{\frac{2}{3}}$ they will be decomposed into series based on whole positive powers of $m$.

Assuming that these decompositions are framed in a series (29), and considering the series obtained in this way, they are arranged by powers of $m, \lambda$ and $e^{2i\tau}$, we'll show that at the values of $m$ and $\lambda$ relevant for lunar theory, and for realistic values of $\tau$, and these series will be absolutely convergent. That is, there will be absolutely converging series, based on powers of $m$, $e^{2i\tau}$, which we see immediately, if we put the power $\lambda = m^2$.

We will thus prove the convergence of the series, based on powers of $m$, which can be represented by coefficients $a_s$.

Finally, if, instead of $m$, we wish to use the value $m_1 = \frac{m}{1+m}$, representing the average movement of the sun to the average stellar movement of the moon, and, just as was done in some old moon theories, the coefficients $a_s$, presented as series based on powers of $m_1$, the convergence of these series will follow at once from the absolute convergence of the series, based on powers of $m$, if we take into account that the decomposition of $m$ by powers of $m_1$, has positive coefficients.

5. Assuming $m$ to be a positive number and considering the formula (27) at different values of $\sigma$ taken from the series

$$\pm 1, \pm 2, \pm 3, ...,$$

we notice at once that the coefficients $A_{r,\sigma}$ and $A_{r,-\sigma}$ reach the highest numerical values in the case of $\sigma = -1$. We're going to therefore have inequality

$$|a_{r,\sigma}| < \frac{\frac{3}{2}l|A_{r,-\sigma}| + \left(8 + 4m + \frac{3}{2}l\right)|A_{r,\sigma}|}{2(6 - 4m + m^2)}, \qquad (34)$$

valid for all the values mentioned now of $\sigma$.

But inequality is also true for $\sigma = 0$.

Indeed, in this case, the second part of it leads to

$$\frac{8 + 4m + 3l}{2(6 - 4m + m^2)}|A_{r,0}|,$$

which represents the large value

$$|a_{r,0}| = \frac{1}{3l}|A_{r,0}|,$$

for all positive values of $m$ the difference

$$\frac{8 + 4m + 3l}{2(6 - 4m + m^2)} - \frac{1}{3l} = \frac{22 + 20m + 9m^2}{4(6 - 4m + m^2)} - \frac{2}{6 + 12m + 9m^2},$$

is obviously positive.

Thus, inequality (34) holds for all values of $\sigma$ that need to be considered here.

Based on this inequality, now it's easy to get a formula giving upper limits of modules of the values $a_{r,\sigma}$, relevant to any $r$, when the upper limits of all modules are known $a_{j,\sigma}$, for which $j < r$.

To do this, we notice that for $A_{r,\sigma}$, there is a polynomial made up of values $a_{j,s}$, where all coefficients are positive. This stems from the very definition of $A_{r,\sigma}$, if we take into account that the coefficients of the decomposition $R_p$ and $R_q$ based on powers of $p$ and $q$ are all positive.

It follows that by replacing in terms of $A_{r,\sigma}$ all $a_{j,s}$ their modules or any upper limits of their modules, we will get some upper limit for $|A_{r,\sigma}|$.

Let's say that in some way we have found upper limits of modules for all $a_{j,s}$ for which $j < r$. Meaning these are upper limits via $a_{j,s}$ and the result of the substitution in terms of $A_{r,\sigma}$ values $a_{j,s}$ values $a_{j,s}$ via $\mathbf{A}_{r,\sigma}$ put

$$\mathbf{a}_{r,\sigma} = \frac{\frac{3}{2}l\mathbf{A}_{r,-\sigma} + \left(8 + 4m + \frac{3}{2}l\right)\mathbf{A}_{r,\sigma}}{2(6 - 4m + m^2)^2}. \tag{35}$$

Then, based on what's been seen now, it's an expression $\mathbf{a}_{r,\sigma}$ will be some upper limit of the module $a_{r,\sigma}$.

So, taking

$$\mathbf{a}_{1,1} = |a_{1,1}|, \qquad \mathbf{a}_{1,-1} = |a_{1,-1}| \tag{36}$$

and identifying everything else $\mathbf{a}_{j,s}$ by the formula (36), we have upper limits of the modules of all $a_{j,s}$.

It's not hard to see that, if you take

$$\mathbf{A}_{1,-1} = \frac{3}{2}, \qquad \mathbf{A}_{1,1} = 0,$$

formula (35) for $r = 1$ will lead to equality (36) and, therefore, it can be used for all values of $r$, that need to be looked at here.

Note that the values determined in this way for $\mathbf{a}_{j,s}$ everything will be increasing functions of $m$, as long as $m$, remaining positive, does not exceed a certain limit. So, for example, this will obviously take place until $m < 2$ for the expression

$$6 - 4m + m^2$$

will be a decreasing function of $m$

Now let's take a look at the series

$$\sum_{r=1}^{r=\infty} \sum_{s=0}^{s=r} \mathbf{a}_{r,r-2s} \lambda^r, \tag{37}$$

and let's say that for any positive values of $m$ and $\lambda$ which we'll call $m_1$ and $\lambda_1$, proven to be convergent. Then if $m_1 < 2$, will be proved the absolute convergence of the series (29) for the values of $\tau$ for which $|\lambda| < \lambda_1$ for all positive values of $m$, which doesn't exceed $m_1$.

But the condition of convergence of the series (37) is not difficult to find, because this series can be defined by some algebraic equation.

Indeed, from (35) follows

$$\sum_{s=0}^{s=r} \mathbf{a}_{r,r-2s} = \frac{8 + 4m + 3l}{2(6 - 4m + m^2)} \sum_{s=0}^{s=r} \mathbf{A}_{r,r-2s},$$

and the sum in the second part of the equality is the result of a substitution in terms of

$$\frac{3}{2}q_{r-1} + lP_r$$

each $p_j$, and each $q_j$ expression

$$\sum_{s=0}^{s=j} \mathbf{a}_{j,j-2s}$$

appropriately, for each $j$.

So, in the series (37) consider each term, a common member of which appears to be an expression

$$\sum_{s=0}^{s=r} a_{r,r-2s}\lambda^r,$$

that series this one will be produced by decomposition by powers of $\lambda$ disappearing at $\lambda = 0$ with root $z$ the equation

$$z = \frac{8 + 4m + 3l}{2(6 - 4m + m^2)}\left\{\frac{3}{2}(1 + z)\lambda + l\left[\frac{1}{(1-z)^2} - 1 - 2z\right]\right\}.$$

Taking into account the expression $l$ and assuming

$$\varepsilon = \frac{3(22 + 20m + 9m^2)}{8(6 - 4m + m^2)}\lambda, \quad h = \frac{22 + 20m + 9m^2}{4(6 - 4m + m^2)}l, \tag{38}$$

the equation will be of the form

$$z = \varepsilon(1 + z) + h\frac{(3 - 2z)z^2}{(1 - z)^2}, \tag{39}$$

and our task will lead thus to finding the conditions of convergence of the series, based on powers of $\varepsilon$, in which, for sufficiently small $\varepsilon$, decomposes the root of the equation (39), disappears at $\varepsilon = 0$.

We have considered so far only positive values of $m$. Let's look at all its meanings now, both physical and complex, whose modules do not exceed a certain limit $M$, and show that with a sufficiently small $M$ formula (35), if in it $m$ is replaced by $M$, will give an upper limit to the modules of the value $a_{r,\sigma}$, for all $m$, module which does not exceed $M$.

We will base it on the next easy-to-prove proposition:

If $a, b, c$ are positive numbers that satisfy the inequality

$$ac - b^2 > 0,$$

and $M$ is a positive number that does not exceed the smaller root of the equation

$$\frac{x^2}{a} - 2\frac{x}{b} + \frac{1}{c} = 0,$$

the lower limit of the module function

$$a - 2bx + cx^2$$

for the value of $x$, modules that don't exceed $M$, there will be a number

$$a - 2bM + cM^2.$$

Applying this statement to the function

$$2(4\sigma^2 - 1) - 4m + m^2,$$

included in the denominator of expression (27), conclude that if $M$ does not exceed the number

$$\frac{2\sqrt{4\sigma^2 - 1}}{\sqrt{4\sigma^2 - 1} + \sqrt{4\sigma^2 - 3}} \tag{40}$$

(where radicals are supposed to be positive) then at $|m| \leq M$ inequality will be implemented

$$|2(4\sigma^2 - 1) - 4m + m^2| \geq 2(4\sigma^2 - 1) - 4M + M^2.$$

Because the number (40) is always larger than one, with an unlimited increase of $\sigma$ approaches a limit equal to 1, then, desiring to ensure that the current inequality is valid for all values of $\sigma$ that have to be considered when dealing with a formula (27), we have to assume $M \leq 1$. In this assumption we will have

$$|2(4\sigma^2 - 1) - 4m + m^2| \geq 6 - 4M + M^2$$

for all values of $m$, modules that do not exceed $M$.

From this we conclude that if the assumption made in formula (35) $m$ is replaced by $M$ (as in coefficients $\mathbf{A}_{r,-\sigma}$, and $\mathbf{A}_{r,\sigma}$, and in the latter's expressions, which contain $l$) and using the labels $\mathbf{a}_{j,s}$, in the expressions $\mathbf{A}_{r,-\sigma}$, and $\mathbf{A}_{r,\sigma}$, we will understand any upper limits of modules of values $a_{j,s}$, for $|m| \leq M$, then at least for $\sigma$ greater than zero this formula will give upper limits to modules of value $a_{r,s}$, on the same condition $|m| \leq M$.

But it's not difficult to see that for sufficiently small $M$ this formula will also give an upper limit for $a_{r,0}$.

Indeed, based on the above proposition, the module function

$$l = 1 + 2m + \frac{3}{2}m^2$$

for $|m| \leq M$ will be no less than

$$1 - 2M + \frac{3}{2}M^2$$

Whenever $M$ does not exceed the value

$$1 - \frac{1}{\sqrt{3}} = 0{,}42 \dots$$

Therefore, assuming

$$M \leq 1 - \frac{1}{\sqrt{3}}, \tag{41}$$

from (28) we have

$$|a_{r,0}| < \frac{2\mathbf{A}_{r,0}}{3(2 - 4M + 3M^2)},$$

so, noticing that the difference

$$\frac{22 + 20M + 9M^2}{4(6 - 4M + M^2)} - \frac{2}{3(2 - 4M + 3M^2)} = \frac{84 - 112M + 4M^2 + 72M^3 + 81M^4}{12(6 - 4M + M^2)(2 - 4M + 3M^2)}$$

is positive for all values of $M$, and defined by formula (35) (after the above replacement) $\mathbf{a}_{r,0}$, we conclude that

$$|a_{r,0}| < \mathbf{a}_{r,0}.$$

Thus, we are convinced that if, on the condition of (41) in the formula (35) $m$ is replaced by $M$ and this formula will be used, starting from $r = 1$, still taking

$$\mathbf{A}_{1,-1} = \frac{3}{2}, \qquad \mathbf{A}_{1,1} = 0,$$

it will give upper limits for modules of value $a_{r,\sigma}$, for all $m$, modules which do not exceed $M$.

Now let's consider the series, based on powers of $m, \lambda$, and $e^{2i\tau}$, which turn the series (29) after the decomposition of the coefficients $a_{r,\sigma}$, by ascending powers of $m$.

Accepting the condition (41) and reasoning under $\mathbf{a}_{r,\sigma}$, the highest limits that have now been discussed, based on the famous theorem, conclude that the decomposition $a_{r,\sigma}$, by ascending powers of $m$ are numerically less than the appropriate ratios of the same decomposition function

$$\mathbf{a}_{r,\sigma}\left(1 - \frac{m}{M}\right)^{-1}.$$

Therefore, to prove the convergence of the series under consideration in the $\tau$ and under the conditions

$$|m| < M, \qquad |\lambda| \leq \lambda_1$$

it will suffice to prove the convergence of the series (37) for $\lambda = \lambda_1$.

Thus, the question of the convergence of these series will lead to the very problem that we have come to, considering the series (29) for $m$ positive.

6. We turn to the equation (39), which we will consider on the assumption that $\varepsilon$ and $h$ are positive numbers.

Our task was to determine the condition of convergence of the series

$$\varepsilon + c_2 \varepsilon^2 + c_3 \varepsilon^3 + \cdots, \tag{42}$$

based on ascending powers of $\varepsilon$, which, with sufficiently small $\varepsilon$, decompose the root of the equation (39), which disappears at $\varepsilon = 0$.

Let's point out first one very simple condition sufficient for this convergence.

We get this condition immediately, if, instead of the equation (39), which is third power relative to $z$, we look at the following square equation:

$$z = \varepsilon(1+z) + \frac{9hz^2}{3-4z}. \tag{43}$$

It's easy to see that the coefficients of the decomposition by ascending powers of $z$ in the second part of this equation are no less than the coefficients of the same decomposition of the second part of the equation (39). Therefore, the series (42), which for sufficiently small $\varepsilon$ will decompose the root of the equation (43), and disappear at $\varepsilon = 0$, will have no smaller coefficients than the series (42), and for the convergence of the latter it will be enough to demonstrate the convergence of this new series.

But the series mentioned here is obtained by decomposition by ascending powers of $\varepsilon$ of the function

$$z = \frac{3 + \varepsilon - \sqrt{9 - 6(7+18h)\varepsilon + 49\varepsilon^2}}{2(4 + 9h - 4\varepsilon)},$$

consideration of which leads to the following condition on its convergence:

$$\varepsilon \leq \frac{3}{7 + 18h + \sqrt{(7+18h)^2 - 49}},$$

where the radical should be considered positive.

From here you can see that if

$$\varepsilon \leq \frac{3}{2(7 + 18h)},$$

and even more so if

$$\varepsilon \leq \frac{1}{6(1+2h)}, \tag{44}$$

the series considered, and therefore, series (42), will converge.

The condition (44) is what we wanted to point out.

Without stopping at the evidence, let's say now, what is the condition, not only sufficient for the convergence of the series (42), but also necessary.

Put

$$\omega = \left(\frac{2h}{1 - \varepsilon + 2h}\right)^{1/3}.$$

Then the condition will be expressed in such a way:

$$\varepsilon \leq \frac{2 - \omega - \omega^2}{4 + \omega + \omega^2} < 1, \tag{45}$$

that obviously requires that at this given $h$ number $\varepsilon$ does not exceed a certain limit $\varepsilon_0$, lying between $0$ and $\frac{1}{2}$.

You can see that the condition is the same as the one that the equation (39) has roots between 0 and 1, and that the series (42) represents the smallest of these roots.

We now turn to expressions $\varepsilon$ and $h$, that are given by formulas (38), and assume
$$\lambda = m^2.$$

At the same time, $\varepsilon$ and $h$ are made by the functions of one $m$, and there is a question of the greatest of $m$ for which the series (42) still converges (and converges for all smaller values of $m$).

Considering the condition (45), we can make sure that the value of this lies between $\frac{1}{7}$ and $\frac{1}{6}$. But the conditions are already (44) enough to show that $\frac{1}{7}$ has not yet reached this limit.

Indeed, when $m = \frac{1}{7}$ formula (38) gives
$$\varepsilon = \frac{1227}{34\,888} = \frac{8589}{244\,216}, \qquad h = \frac{53\,761}{34\,888},$$

from which follows
$$\frac{1}{6(1+2h)} = \frac{8\,722}{210\,615} > \varepsilon.$$

Thus, we are sure that when $m \leq \frac{1}{7}$ all series we reviewed will converge.

For the theory of the moon this conclusion is more than enough, because the value of $m$, corresponding to the moon, does not exceed $\frac{1}{12}$.

7. The equation (39), which served us to solve the issue of convergence of series (29) and (30), can also serve to determine the higher margin of error, that we're doing, discarding in these series all members any power of $\lambda$, for discarding in any of these series all members above the $n$th power, we get the result, the margin of error, obviously doesn't exceed the value
$$z - (\varepsilon + c_2\varepsilon^2 + \cdots + c_n\varepsilon^n), \tag{46}$$
where $z$ is the smallest root of the equation (39).

It cannot, however, be expected to have the highest degree of accuracy that is desirable in the task at hand.

The reasons for its lack of accuracy are, of course, at the very heart of our analysis, and the most important ones are hardly easily eliminated. However, one of the reasons for the inaccuracy can be easily addressed.

In fact, to the equation (39) we came across by considering a series (37) in the assumption that
$$\mathbf{a}_{1,-1} = |a_{1,-1}|, \qquad \mathbf{a}_{1,1} = |a_{1,1}|$$
and that all other $\mathbf{a}_{r,\sigma}$, starting from $r = 2$, determined by formula (35). But we could only use this last formula for $r > N$, reasoning under $N$ some number, greater than 1, and take
$$\mathbf{a}_{r,\sigma} = |a_{r,\sigma}|$$
for $r \leq N$. Then the rest of the series (37) after members of the $n$th power of $\lambda$ would obviously present a higher margin of error than the above value of (46).

Let's consider the same series (37) in this new assumption.

It is easy to see that this series can still be defined by some third-degree algebraic equation.

To make up this equation, let's say that in a function
$$\varepsilon z + h \frac{(3-2z)z^2}{(1-z)^2}$$
for the values (38) for $\varepsilon$ and $h$ we use instead of $z$ the expression
$$\sum_{j=1}^{j=N} \sum_{s=0}^{s=j} |a_{j,j-2s}| \lambda^j$$

and taking into account the dependency of $\varepsilon$ on $\lambda$, the result is decomposed into a series

$$L_2\lambda^2 + L_3\lambda^3 + \cdots + L_N\lambda^N + \cdots$$

on ascending posers of $\lambda$. Then, assuming

$$\varepsilon' = \sum_{j=1}^{j=N}\sum_{s=0}^{s=j}|a_{j,j-2s}|\lambda^j - \sum_{j=2}^{j=N}L_j\lambda^j,$$

we can form the equation like this:

$$z = \varepsilon' + \varepsilon z + h\frac{(3-2z)z^2}{(1-z)^2}. \tag{47}$$

Let's say now that

$$l_1\lambda + l_2\lambda^2 + l_3\lambda^3 + \cdots \tag{48}$$

there's a series which decomposes by powers of $\lambda$ and the root of this equation disappears for $\lambda = 0$. Then the residual error, which we get, stopping, when calculating any of the series (30), at members $n$th power of $\lambda$, won't exceed the value

$$z - (l_1\lambda + l_2\lambda^2 + \cdots + l_n\lambda^n),$$

where $z$ the smallest root of the equation (47).

This way we get the largest margin of error, which when $N$ is sufficiently large can achieve proper accuracy. But it's possible that for this we must take $N > n$, as is the case in the first example cited below.

To obtain the higher limit in question, it is necessary to calculate the lowest root of the equation (47), for which we can resort to the method of successive approximations, bringing the pre-equation to a simpler form

$$z = \delta + g\frac{(3-2z)z^2}{(1-z)^2}, \tag{49}$$

where

$$\delta = \frac{\varepsilon'}{1-\varepsilon}, \quad g = \frac{h}{1-\varepsilon}.$$

The assumption, implied everywhere, is that $\lambda, \varepsilon$ and $h$ are positive numbers, satisfying the convergence of the series (48), we can prove that $\delta$ and $g$ will be positive numbers, satisfying the inequality

$$\frac{1+2\delta}{3(1+2g)} < 1,$$

and that the smallest root of the equation (49) (which when the assumption is made, have only physical and positive roots) will not exceed the number

$$\frac{1+2\delta}{3(1+2g)}, \tag{50}$$

not surpassing either of the two roots.

Under these conditions, if $z_0$ is some positive number, not superior to (50), and by that number $z_0$ consistently determines $z_1, z_2, z_3$ etc., taking advantage of the equation

$$z_{n+1} = \delta + g\frac{(3-2z_n)z_n^2}{(1-z_n)^2},$$

that $z_n$, for increasing $n$, will approach the root of the $c$, or constantly increasing (when $z_0 < c$), or constantly decreasing (when $z_0 > c$).

In this way, the root we are interested in can be calculated fairly quickly with great accuracy.

8. I will now apply the theory outlined to lunar theory, limiting itself to the assumption that $N = 2$. In this assumption, noting that

$$L_2\lambda^2 = c_2\varepsilon^2 = (1 + 3h)\varepsilon^2,$$

and assuming

$$\varepsilon_1 = \{|a_{2,0}| + |a_{2,2}| + |a_{2,-2}|\}\lambda^2,$$

we find

$$\varepsilon' = \varepsilon + \varepsilon_1 - (1 + 3h)\varepsilon^2,$$
$$l_1\lambda = \varepsilon, \quad l_2\lambda^2 = \varepsilon_1, \quad l_3\lambda^3 = (1 + 6h)\varepsilon\varepsilon_1 + 4h\varepsilon^3,$$

To calculate these expressions, you need to find $a_{2,0}, a_{2,-2}$ and $a_{2,2}$.

When known $a_{1,-1}$ and $a_{1,1}$ that will be found by formulas

$$a_{1,-1} = \frac{3}{16}\frac{38 + 28m + 9m^2}{6 - 4m + m^2}, \quad a_{1,1} = -\frac{9}{16}\frac{2 + 4m + 3m^2}{6 - 4m + m^2},$$

and calculated

$$2l = 2 + 4m + 3m^2,$$

$a_{2,0}$ found by formula

$$a_{2,0} = -\frac{1}{4}(a_{1,-1} - a_{1,1})^2 - \left(2a_{1,1} + \frac{1}{2l}\right)a_{1,-1},$$

which is easily seen from (28), and from (31) and (32).

As for $a_{2,-2}$ and $a_{2,2}$ then to calculate them, instead of formula (27), we turn to the formulas of Hill, which are more likely to lead to the goal.

From (31) and (33) we find

$$a_{2,-2} = -\alpha_{-2,0}, \quad a_{2,2} = -\alpha_{2,0},$$

and from the formulas in part 2

$$\alpha_{2,0} = 2\overline{[2]}\alpha_{1,0} + [2,1]\alpha_{1,0}\alpha_{-1,0},$$
$$\alpha_{-2,0} = 2\overline{(-2)}\alpha_{1,0} + [-2,-1]\alpha_{1,0}\alpha_{-1,0}.$$

Hence, taking into account that

$$\alpha_{-1,0} = -a_{1,-1}, \quad \alpha_{1,0} = -a_{1,1},$$
$$\overline{[2]} = \frac{3}{64}\frac{2 + 16m + 9m^2}{30 - 4m + m^2}, \quad \overline{(-2)} = -\frac{9}{64}\frac{38 + 16m + 3m^2}{30 - 4m + m^2},$$
$$[2,1] = -\frac{1}{2}\frac{26 + m^2}{30 - 4m + m^2}, \quad [-2,-1] = -\frac{1}{2}\frac{18 - 8m + m^2}{30 - 4m + m^2}.$$

we have

$$a_{2,-2} = \frac{16(18 - 8m + m^2)a_{1,-1} - 9(38 + 16m + 3m^2)a_{1,1}}{32(30 - 4m + m^2)},$$

$$a_{2,2} = \frac{16(26 + m^2)a_{1,-1} - 3(2 + 16m + 9m^2)a_{1,1}}{32(30 - 4m + m^2)}.$$

According to these formulas, taking $\lambda = m^2$, we have

$$a_{1,-1}\lambda, \quad a_{1,1}\lambda, \quad a_{2,-2}\lambda^2, \quad a_{2,2}\lambda^2,$$

that will deliver immediately

$$\varepsilon = \{|a_{1,-1}| + |a_{1,1}|\}\lambda$$

and $\varepsilon_1$, and then calculating $h$ from formula (38), we have $l_3\lambda^3, \varepsilon', \delta$ and $g$.

We'll take with Hill.

$$m = 0.08084\ 89338\ 08312$$

and when you write out each calculation result, we will keep the usual requirement, so that the error does not exceed half of the unit of the last decimal place at which we stop. Limited to ten decimal places, we have

$$\begin{aligned}
a_{1,-1}\lambda &= \phantom{-}0.00869\ 58085\ (-)^4, & h &= 1.22011\ 19633\ (-), \\
a_{1,1}\lambda &= -0.00151\ 58492\ (-), & l_1\lambda &= 0.01021\ 16577\ (1), \\
\varepsilon &= \phantom{-}0.01021\ 16577\ (-), & l_2\lambda^2 &= 0.00003\ 00097\ (+), \\
a_{2,-2}\lambda^2 &= -0.00000\ 01637\ (-), & l_3\lambda^3 &= 0.00000\ 77466\ (-), \\
a_{2,0}\lambda^2 &= -0.00002\ 39667\ (+), & \varepsilon' &= 0.00975\ 60241\ (-), \\
a_{2,2}\lambda^2 &= -0.00000\ 58793\ (+), & \delta &= 0.00985\ 66771\ (-), \\
\varepsilon_1 &= \phantom{-}0.00003\ 00097\ (+), & g &= 1.23269\ 98742\ (-).
\end{aligned}$$

Turning now to the equation (49) and limited to calculating the smallest root to eight decimal places, we find

$$z = 0.01025064\ (-).$$

From here, with the same approximation

$$z - l_1\lambda = 0.00003\ 898, \tag{51}$$

$$z - l_1\lambda - l_2\lambda^2 = 0.00000897, \tag{52}$$

$$z - l_1\lambda - l_2\lambda^2 - l_3\lambda^3 = 0.00000\ 12^{3(-)}_{2(+)}. \tag{53}$$

The numbers represent the highest margin of error, related to the cases $n = 1, n = 2$ and $n = 3$.

Comparing the number (51) with value $\varepsilon_1$, which can be achieved by the second-order members in the expression of function $\xi$ [formula (30)], we conclude that the number is quite satisfactory in accuracy by the highest margin of error.

To judge the accuracy of the upper limit (52), we need to calculate length of expressions (30) at lease the coefficients of order relative to $\lambda$, and make yourself some idea of the size

$$\varepsilon_2 = \{|a_{3,-3} - a_{1,-1}| + |a_{1,1}| + |a_{3,3}|\}\lambda^3.$$

Turning to formulas (31) and (33), we have

$$a_{3,\pm 3} = -\alpha_{\pm 3,0}, \quad a_{3,\pm 1} = -a_{2,0}\alpha_{1,\pm 1} - \alpha_{\pm 1,1},$$

We must then calculate

$$a_{2,0}a_{1,\pm 1}\lambda^3, \quad \alpha_{\pm 1,1}\lambda^3, \quad \text{and} \quad \alpha_{\pm 3,0}\lambda^3.$$

From the above numerical results, we deduce

$$\begin{aligned}
a_{2,0}a_{1,-1}\lambda^3 &= -0.00000\ 02085\ (-), \\
a_{2,0}a_{1,1}\lambda^3 &= \phantom{-}0.00000\ 00414\ (-).
\end{aligned}$$

To get the same $\alpha_{1,\pm 1}\lambda^3$ and $\alpha_{\pm 3,0}\lambda^3$ turn to Hill's numerical findings, from which we borrow

$$\begin{aligned}
\alpha_{-1,1}\lambda^3 &= \phantom{-}0.00000\ 00616\ (-), \\
\alpha_{1,1}\lambda^3 &= -0.00000\ 01417\ (-),
\end{aligned}$$

---

[4] The sign $(-)$ indicates that the exact value of the number is less than write, and the sign $(+)$, that it is more than written.

$$\begin{aligned} \alpha_{-3,0}\lambda^3 &= 0.00000\ 00025\ (-), \\ \alpha_{3,0}\lambda^3 &= 0.00000\ 00300\ (+) \end{aligned}$$

From here we find

$$\begin{aligned} a_{3,-1}\lambda^3 &= 0.00000\ 01469, \\ a_{3,1}\lambda^3 &= 0.00000\ 01003, \\ a_{3,-3}\lambda^3 &= -0.00000\ 00025\ (-), \\ a_{3,3}\lambda^3 &= -0.00000\ 00300\ (+) \end{aligned}$$

and finally,

$$\varepsilon_2 = 0.00000\ 0280\ (-).$$

So

$$\varepsilon_2 < 0.00000028,$$

And the number (53) does not exceed

$$0.00000\ 123.$$

Therefore, the totality of members above the second order in the series (30) probably will be numerically less than

$$0.00000\ 151, \qquad (54)$$

and therefore, our upper limit (52), almost 6 times higher than that number, cannot be called satisfactory.

Even less satisfactory will be, of course, the upper limit (53).

In conclusion, we will use the numbers to show the accuracy of some of the approximate expressions $\xi$ and $\eta$.

Let's say that in the series (30) we neglected all members above the first order regarding $\lambda$ and that the coefficients members of the first order are calculated to five decimal places, which gives

$$\xi = -0.00718\cos 2\tau, \quad \eta = 0.01021\sin 2\tau.$$

Taking the highest margin of error for the case $n = 2$ the number (54) and combining it with $\varepsilon_1$, get for the upper margin of the error in the case $n = 1$ the following number:

$$0.00003152.$$

Therefore, noting that the errors of the coefficients of our approximate expressions $\xi, \eta$, do not reach 0.000002, and conclude that the calculation of these expressions, for all values of $\tau$, can be found for $\xi$ and $\eta$ values, whose errors will not exceed 0.00004.

Wanting to achieve more accuracy for the coefficients, calculated to five decimal places, we have to take into account members of the second order. At the same time, stopping at formulas

$$\begin{aligned} \xi &= 0.00002 & \eta &= 0.01021\sin 2\tau \\ &\phantom{=}-0.00718\cos 2\tau & &\phantom{=}+0.00001\sin 4\tau, \\ &\phantom{=}+0.00001\cos 4\tau \end{aligned}$$

we can expect to get, for $\xi$ and $\eta$, values for which the margin of error will not exceed 0.00001.

Ляпунов, А. (1896). О рядах, предложенных Хиллом для представления движения Луны. Труды Отделения физических наук Общества любителей естествознания, М., 1896, т. VIII, вып. 1, 1—23.

Thomas S. Ligon, http://orcid.org/0000-0002-4067-876X

Last change: 2020-11-09.

This is the original Russian version, except that the document was created by traditional typesetting and later scanned to create the PDF. Then I went through the following steps:
- Opened the PDF document in Microsoft Word, which automatically ran OCR (optical character recognition), and converted it to an editable Word document.
- Recreated all formulas and equations using the built-in Microsoft Equation Editor.

О РЯДАХ, ПРЕДЛОЖЕННЫХ ХИЛЛОМ ДЛЯ ПРЕДСТАВЛЕНИЯ
ДВИЖЕНИЯ ЛУНЫ

1. В первом томе «American Journal of Mathematics» в 1878 Г. Хилл опубликовал весьма интересный мемуар «Researches in the Lunar Theory»[5] где задача п движении Луны под действием Земли и Солнца трактуется как предельный случай задачи трех тел, к которому можно перейти, предполагая массу Солнца и его расстояние от Земли беспредельно возрастающими, притом так, чтобы отношение куба этого расстояния к массе Солнца приближалось к некоторому пределу, не зависящему от времени. Решение задачи в этой постановке, на которую, имен в виду дальнейшее исследование вопроса, Хилл смотрит только как на предварительную, доставляет такое приближение, в котором принимается в расчет среднее движение Солнца, но пренебрегается его параллаксом и эксцентриситетом. При большей части своих выводов Хилл пренебрегает также наклонностью лунной орбиты и, предполагая таким образом, что Луна движется плоскости эклиптики, для своих дифференциальных уравнений, которыми он определяет относительное движение Луны по отношению к осям, вращающимся вокруг Земли вместе с Солнцем, ищет не общий интеграл, а некоторое частное решение, подчиненное условию периодичности. Это частное решение, содержащее два произвольных постоянных — один из моментов соединений и среднее движение Луны, — при дальнейших изысканиях Хилла должно было играть роль первого приближения, и оно-то главным образом и изучается в рассматриваемом мемуаре[6]

Свое частное решение Хилл определяет тригонометрическими рядами, в которых коэффициенты, известным образом зависящие от средних движений Солнца и Луны, вычисляются им путем последовательных приближений. Но, занимаясь вопросом о вычислении этих коэффициентов, Хилл вовсе не касается вопроса о сходимости своих рядов.

---

[5] Весьма обстоятельный анализ этого мемуара можно найти в третьем томе сочинения: Tisserand, Traité Mécanique Céleste.

[6] Недавно в «Astronomical Journal» (vol. XV, No. 353) появился новый мемуар Хилла о движении Луны, в котором Хилл рассматривает подобное же периодическое решение, принимая в расчет параллакс Солнца.

Некоторые общие соображения приводят к заключению, что ряды, о которых идет речь, несомненно должны быть сходящимися, если отношение среднего движения Солнца к среднему движению Луны не превосходит некоторого предела, при каковом условии рядами этими представляется одно из периодических движений, существование которых было обнаружено Пуанкаре при его общих изысканиях о периодических решениях дифференциальных уравнений задачи трех тел. Но названный сейчас предел, сколько мне известно, никем не разыскивался, и не было еще доказано, что он достаточно велик для сходимости рассматриваемых рядов при величинах средних движений Солнца и Луны, принимаемых в астрономии.

Между тем в виду важного значения, которое могут иметь ряды Хилла в теории Луны, вопрос о сходимости их заслуживает серьезного внимания. Я счел поэтому небесполезным опубликовать свое исследование вопроса, которое и предлагаю в настоящей статье. Исследование это показывает, что величина отношения средних движений Солнца и Луны, принимаемая в астрономии, вполне обеспечивает сходимость рядов Хилла.

2. Рассматриваемые Хиллом дифференциальные уравнения движения Луны суть следующие:

$$\left.\begin{aligned}\frac{d^2x}{dt^2} - 2n\frac{dy}{dt} + \frac{\mu x}{(x^2+y^2)^{\frac{3}{2}}} &= 3n^2x, \\ \frac{d^2y}{dt^2} + 2n\frac{dx}{dt} + \frac{\mu y}{(x^2+y^2)^{\frac{3}{2}}} &= 0.\end{aligned}\right\} \quad (1)$$

Здесь $x$ и $y$ суть прямоугольные координаты Луны по отношению к осям, начало которых совпадает с центром Земли и которые движутся в плоскости эклиптики так, что ось x постоянно направлена от Земли к Солнцу, а ось y — в сторону годичного движения Солнца относительно Земли; n означает среднее движение Солнца и $\mu$ — сумму масс Земли и Луны, умноженную на коэффициент притяжения.

Этим уравнениям Хилл старается удовлетворит рудами вида

$$x = A_0 \cos n_1(t-t_0) + A_1 \cos 3n_1(t-t_0) + A_2 \cos 5n_1(t-t_0) + \cdots$$
$$y = B_0 \sin n_1(t-t_0) + B_1 \sin 3n_1(t-t_0) + B_2 \sin 5n_1(t-t_0) + \cdots$$

где $t_0$ и $n_1$, — произвольные постоянные, а $A_0, A_1$ и т. д. $B_0, B_1$ и т. д. — постоянные, определяемые, как функции $n$ и $n_1$, из условия задачи.

Постоянное $n_1$ должно представлять синодическое среднее движение Луны.

Полагая

$$n_1(t-t_0) = \tau$$

и делая вообще для всякого $s$

$$A_s = a_s + a_{-s-1}, \quad B_s = a_s - a_{-s-1},$$

Хилл представляет свои ряды под видом

$$x = \sum a_s \cos(2s+1)\tau, \quad y = \sum a_s \sin(2s+1)\tau,$$

где суммирование распространяется на все целые значения $s$ от $-\infty$ до $+\infty$.

Вместо переменных $x$ и $y$ Хилл вводит затем следующие:

$$u = x + y\sqrt{-1}, \quad v = x - y\sqrt{-1},$$

для которых предыдущие ряды дают

$$u = \sum a_s\, e^{(2s+1)i\tau}, \quad v = \sum a_{-s-1}\, e^{(2s+1)i\tau}, \tag{2}$$

где $i = \sqrt{-1}$, а дифференциальные уравнения, если положим

$$\frac{n}{n_1} = m, \quad \frac{\mu}{n_1^2} = k,$$

принимают следующий вид:

$$\left.\begin{aligned}
\frac{d^2 u}{d\tau^2} + 2mi\frac{du}{d\tau} + \frac{ku}{(uv)^{\frac{3}{2}}} &= \frac{3}{2}m^2(u+v), \\
\frac{d^2 v}{d\tau^2} - 2mi\frac{dv}{d\tau} + \frac{kv}{(uv)^{\frac{3}{2}}} &= \frac{3}{2}m^2(u+v).
\end{aligned}\right\} \tag{3}$$

Означая через $C$ произвольное постоянное, из уравнений этих выводим следующее:

$$\frac{du}{d\tau}\frac{dv}{d\tau} - \frac{2k}{\sqrt{uv}} = +\frac{3}{4}m^2(u+v)^2 - C, \tag{4}$$

представляющее уравнение живой силы для рассматриваемой задачи.

Исключая из уравнений (3) и (4) постоянное $k$, находим

$$\left.\begin{aligned}
v\frac{d^2 u}{d\tau^2} - u\frac{d^2 v}{d\tau^2} + 2mi\frac{d(uv)}{d\tau} + \frac{3}{2}m^2(u^2 - v^2) &= 0, \\
v\frac{d^2 u}{d\tau^2} + u\frac{d^2 v}{d\tau^2} + \frac{du}{d\tau}\frac{dv}{d\tau} + 2mi\left(v\frac{du}{d\tau} - u\frac{dv}{d\tau}\right) - \frac{9}{4}m^2(u+v)^2 + C &= 0,
\end{aligned}\right\} \tag{5}$$

Чтобы удовлетворить уравнениям (3), Хилл старается сначала удовлетворить уравнениям (5), и так как всякие функции $u$ и $v$, удовлетворяющие этим последним уравнениям, необходимо удовлетворяют следующим:

$$\frac{d^2 u}{d\tau^2} + 2mi\frac{du}{d\tau} + \frac{hu}{(uv)^{\frac{3}{2}}} = \frac{3}{2}m^2(u+v),$$

$$\frac{d^2 v}{d\tau^2} - 2mi\frac{dv}{d\tau} + \frac{hv}{(uv)^{\frac{3}{2}}} = \frac{3}{2}m^2(u+v),$$

где $h$ — некоторое постоянное, то для окончательного решения задачи ему остается только удовлетворить какому-либо соотношению между постоянными, приводящему к равенству

$$h = k. \tag{6}$$

Разыскивая условия, при которых ряды (2) удовлетворяют уравнениям (5), мы получаем выражение постоянного $C$ через коэффициенты $a_s$, и следующие соотношения между последними:

$$4j \sum_s (2s - j + 1 + m) a_s a_{s-j} - \frac{3}{2} m^2 \sum_s a_s (a_{j-s-1} - a_{-j-s-1}) = 0,$$

$$\sum_s \left[ (2s+1)(2s - 2j + 1) + 4j^2 + 4(2s - j + 1)m + \frac{9}{2} m^2 \right] a_s a_{s-j}$$
$$+ \frac{9}{4} m^2 \sum_s a_s (a_{j-s-1} - a_{-j-s-1}) = 0,$$

приводящиеся к виду

$$\sum_s [8s^2 - 8(4j - 1)s + 20j^2 - 16j + 2 + 4(4s - 5j + 2)m + 9m^2] a_s a_{s-j} + 9m^2 \sum_s a_s a_{j-s-1} = 0,$$

$$\sum_s [8s^2 + 8(2j + 1)s - 4j^2 + 8j + 2 + 4(4s + j + 2)m + 9m^2] a_s a_{s-j} + 9m^2 \sum_s a_s a_{-j-s-1} = 0,$$

Здесь суммирование распространяется на все целые значения $s$ от $-\infty$ до $+\infty$, и $j$ можно давать всякие целые значения, за исключением нуля. Но давая как положительные, так и отрицательные значения, можно ограничиться рассмотрением только одного из двух написанных сейчас уравнений, ибо нетрудно видеть, что одно выводится из другого заменою $j$ на $-j$.

Хилл рассматривает некоторую комбинацию этих уравнений.

Если первое из них представим равенством $L_j = 0$ и, следовательно, второе — равенством $L_{-j} = 0$, то комбинацию эту можем представить так:

$$[4j^2 - 8j - 2 - 4(j + 2)m - 9m^2] L_j + [20j^2 - 16j + 2 - 4(5j - 2)m + 9m^2] L_{-j} = 0.$$

Получаемое таким путем уравнение легко приводится к виду

$$\sum_s ([j, s] a_s a_{s-j} + [j] a_s a_{j-s-1} + (j) a_s a_{-j-s-1}) = 0, \qquad (7)$$

где, по обозначению Хилла,

$$[j, s] = -\frac{s}{j} \frac{4s(j - 1) + 4j^2 + 4j - 2 - 4(s - j + 1)m + m^2}{2(4j^2 - 1) - 4m + m^2},$$

$$[j] = -\frac{3m^2}{16j^2} \frac{4j^2 - 8j - 2 - 4(j + 2)m - 9m^2}{2(4j^2 - 1) - 4m + m^2},$$

$$(j) = -\frac{3m^2}{16j^2} \frac{20j^2 - 16j + 2 - 4(5j - 2)m + 9m^2}{2(4j^2 - 1) - 4m + m^2}.$$

Уравнением (7) Хилл пользуется для определения отношений $\frac{a_s}{a_0}$ ($s = \pm 1, \pm 2, \pm 3, \ldots$), допуская, что отношения эти как функции параметра $m$ уничтожаются вместе с $m$ и что при $m$ бесконечно малом первого порядка $\frac{a_s}{a_0}$ есть бесконечно малая порядка $\pm 2s$. Допущение это позволяет определять рассматриваемые отношения путем последовательных приближений с точностью до членов любого порядка относительно $m$, что окончательно приводит к представлению этих отношений под видом бесконечных рядов.

Укажем здесь ряды, к которым можно таким образом прийти из рассмотрения способа последовательных приближений, употребляемого Хиллом.

Вычисления Хилла основаны на том обстоятельстве, что в уравнении (7) коэффициенты $[j]$ и $(j)$ уничтожаются при $m = 0$ и суть второго порядка относительно $m$, а коэффициенты $[j, s]$ таковы, что

$$[j, j] = -1, \quad [j, 0] = 0. \tag{8}$$

Выделяя множитель $m^2$, мы положим

$$[j] = \overline{[j]}m^2, (j) = \overline{(j)}m^2$$

и затем, вводя параметр $\lambda$, уравнение (7) заменим следующим:

$$\sum_s \left([j, s]a_s a_{s-j} + \overline{[j]}\lambda a_s a_{j-s-1} + \overline{(j)}\lambda a_s a_{-j-s-1}\right) = 0. \tag{9}$$

Полагая теперь

$$\frac{a_s}{a_0} = \lambda^{|s|}\alpha_s, \tag{10}$$

где $|s|$ означает численное значение $s$, будем искать величины $\alpha_s$ ($s = \pm 1, \pm 2, \pm 3, \ldots$) под видом рядов, расположенных по целым положительным степеням $\lambda$.

Подставляя эти ряды в выражение (10) и заменяя затем $\lambda$ через $m^2$, и получим те ряды, к которым приводит способ вычисления Хилла, если при составлении каждого приближения удерживать лишь члены, порядки которых ниже порядка соответствующей погрешности.

Нетрудно видеть, что ряды, которыми представятся величины $\alpha_s$, будут содержать только четные степени $\lambda$ и что коэффициенты в этих рядах будут определяться в известной последовательности так, что для определения каждого коэффициента в каждом ряду достаточно будет конечного числа алгебраических действий.

В самом деле, после положения (10) в уравнении (9) будут встречаться только те степени $\lambda$, показатели которых заключаются в формах

$$|s| + |s - j|, \quad |s| + |s - j + 1| + 1, \quad |s| + |s + j + 1| + 1,$$

а последние, очевидно, дают только числа, равные $|j|$ или превосходящие $|j|$ на числа четные.

Поэтому, если после подстановки (10) все члены уравнения (9) разделить на $\lambda^{|j|}$ то в нем будут встречаться только положительные и четные степени $\lambda$, и вследствие (8) уравнению этому можно будет дать следующий вид:

$$\alpha_j = L_{j,0} + \sum_{s=1}^{\infty} L_{j,s} \lambda^{2s},$$

где

$$L_{1,0} = \overline{[1]}, \quad L_{-1,0} = \overline{(-1)}$$

и вообще

$$L_{j,0} = \begin{cases} \overline{[j]} \sum_{s=0}^{j-1} \alpha_s \alpha_{j-s-1} + [j,s] \sum_{s=1}^{j-1} \alpha_s \alpha_{s-j}, \text{если } j > 1, \\ \overline{(j)} \sum_{s=0}^{-j-1} \alpha_s \alpha_{-j-s-1} + [j,s] \sum_{s=j+1}^{-1} \alpha_s \alpha_{s-j}, \text{если } j < -1, \end{cases}$$

притом всякое $L_{j,s}$ есть совокупность членов вида

$$C \alpha_p \alpha_q,$$

для которых значки $p$ и $q$ удовлетворяют равенствам

$$p - q = j, \quad |p| + |q| = 2s + |j|$$

или

$$p + q = \pm j - 1, \quad |p| + |q| = 2s + |j| - 1.$$

Отсюда ясно, что рассматриваемые ряды будут вида

$$\alpha_s = \alpha_{s,0} + \alpha_{s,1} \lambda^2 + \alpha_{s,2} \lambda^4 + \cdots$$

и что коэффициенты $\alpha_{s,j}$ можно расположить в такой линейный ряд, в котором каждый член будет определяться по предшествующим, причем первым членом будет какой-либо из коэффициентов $\alpha_{1,0}, \alpha_{-1,0}$ определяемых непосредственно по формулам

$$\alpha_{1,0} = \overline{[1]}, \quad \alpha_{-1,0} = \overline{(-1)}.$$

Можно заметить, что каждый из коэффициентов $\alpha_{s,j}$ будет рациональной функцией $m$, в которой, если ее представить под видом несократимой дроби, знаменатель будет содержать только множители типа

$$2(4\sigma^2 - 1) - 4m + m^2, (\sigma = 1,2,3, \dots),$$

не обращающиеся в нуль ни при каких вещественных значениях $m$

Определив отношение $\frac{a_s}{a_0}$, Хилл переходит к разысканию $a_0$ и, составляя различными способами уравнение (6), получает несколько различных выражений для $a_0$. Одно из этих выражений, выводимое из рассмотрения членов с первой степенью $e^{i\tau}$ в результате подстановки рядов (2) в первое из уравнений (3), есть следующее:

$$a_0 = \frac{(kJ)^{1/3}}{\left(1 + 2m + \frac{3}{2}m^2 + \frac{3}{2}m^2 \frac{a_{-1}}{a_0}\right)^{1/3}}, \tag{11}$$

где $J$ означает не зависящий от $\zeta$ член в разложении выражения

$$\left\{\sum \frac{a_s}{a_0} \zeta^{2s}\right\}^{-1/2} \left\{\sum \frac{a_s}{a_0} \zeta^{-2s}\right\}^{-3/2}$$

по положительным и отрицательным степеням $\zeta$.

Применяя свои формулы к теории Луны, Хилл принимает

$$m = 0{,}08084\,89338\,08312$$

и при окончательных выводах пренебрегает лишь членами выше 13-го порядка относительно $m$. Вычисляя при этом отношения $\frac{a_s}{a_0}$ с пятнадцатью десятичными знаками, он считает погрешности не превосходящими двух единиц последней десятичной. Это заключение, впрочем, основано только на сопоставлении вычисленных членов различных порядков и не может считаться доказанным, так как на самом деле Хилл не дает никаких средств для определения высших пределов погрешностей при вычислении отношений $\frac{a_s}{a_0}$ по предлагаемому им способу.

Изложивши сущность анализа Хилла, переходим к вопросу о сходимости его рядов, которые будем определять, однако, несколько иначе.

3. Нетрудно убедиться, что уравнение (9), которым, введя параметр $\lambda$, мы заменили выше уравнение (7), получается при разыскании рядов вида (2), удовлетворяющих следующим дифференциальным уравнениям:

$$\left.\begin{aligned}\frac{d^2 u}{d\tau^2} + 2mi\frac{du}{d\tau} - \frac{3}{2}m^2 u + \frac{ku}{(uv)^{\frac{3}{2}}} &= \frac{3}{2}\lambda v,\\ \frac{d^2 v}{d\tau^2} - 2mi\frac{dv}{d\tau} - \frac{3}{2}m^2 v + \frac{kv}{(uv)^{\frac{3}{2}}} &= \frac{3}{2}\lambda u,\end{aligned}\right\} \tag{12}$$

которые, по преобразовании к переменным $x$ и $y$, принимают вид

$$\left.\begin{aligned}\frac{d^2 x}{d\tau^2} - 2m\frac{dy}{d\tau} - \frac{3}{2}m^2 x + \frac{kx}{(x^2+y^2)^{\frac{3}{2}}} &= \frac{3}{2}\lambda x,\\ \frac{d^2 y}{d\tau^2} + 2m\frac{dx}{d\tau} - \frac{3}{2}m^2 y + \frac{ky}{(x^2+y^2)^{\frac{3}{2}}} &= -\frac{3}{2}\lambda y.\end{aligned}\right\} \tag{13}$$

Эти уравнения мы теперь и будем рассматривать вместо уравнений (1), причем будем иметь в виду, что в окончательном результате должно полагать $\lambda = m^2$.

Делая в уравнениях (13) $\lambda = 0$, получим уравнения, которым можно удовлетворить, принимая

$$x = a\cos\tau, y = a\sin\tau \tag{14}$$

и разумея под $a$ постоянное, определяемое равенством

$$\frac{k}{a^3} = 1 + 2m + \frac{3}{2}m^2.$$

Это постоянное представит длину, мало отличающуюся от среднего расстояния Луны от Земли, куб которого равен

$$\frac{k}{(1+m)^2}.$$

Возвращаясь теперь к уравнениям (13), в которых $\lambda$ будем предполагать численно не превосходящим некоторого предела, будем для них искать решение, в котором $x$ и $y$ были бы периодическими функциями $\tau$ с периодом $2\pi$, обращающимися при $\lambda = 0$ в выражения (14).

Мы будем таким образом искать периодическое движение, которое при малых значениях $\lambda$ мало отличалось бы от равномерного движения по кругу радиуса $a$ с центром в центре Земли.

Искомое решение мы можем представить формулами

$$x = a(1+\xi)\cos\tau - a\eta\sin\tau, \quad y = a(1+\xi)\sin\tau + a\eta\cos\tau,$$

разумея под $\xi$ и $\eta$, периодическое функции $\tau$, уничтожающиеся при $\lambda = 0$.

Вместо этих функций мы будем далее рассматривать следующие их комбинации:

$$-\xi - i\eta = p, \quad -\xi + i\eta = q,$$

при посредстве которых введенные выше переменные $u$ и $v$ выразятся формулами

$$u = a(1-p)e^{i\tau}, \quad v = a(1-q)e^{-i\tau},$$

а в силу последних уравнения (12) доставят следующие дифференциальные уравнения для определения функций $p$ и $q$:

$$\frac{d^2p}{d\tau^2} + 2(1+m)i\frac{dp}{d\tau} + l(1-p) - \frac{l}{(1-p)^{1/2}(1-q)^{3/2}} = \frac{3}{2}\lambda e^{-2i\tau},$$

$$\frac{d^2q}{d\tau^2} - 2(1+m)i\frac{dq}{d\tau} + l(1-q) - \frac{l}{(1-p)^{3/2}(1-q)^{1/2}} = \frac{3}{2}\lambda e^{+2i\tau},$$

где положено для сокращения

$$l = 1 + 2m + \frac{3}{2}m^2.$$

Входящие сюда функции

$$(1-p)^{-1/2}(1-q)^{-3/2}, \quad (1-p)^{-3/2}(1-q)^{-1/2}$$

мы будем представлять рядами, расположенными по целым положительным степеням $p$ и $q$, так что анализ наш будет относиться только к случаю, когда модули функций $p$ и $q$ остаются всегда меньшими единицы.

Означая в названных рядах совокупности членов выше первого измерения относительно $p$ и $q$ через $R_p$ и $R_q$, будем иметь

$$(1-p)^{-1/2}(1-q)^{-3/2} = 1 + \frac{1}{2}p + \frac{3}{2}q + R_p$$
$$(1-p)^{-3/2}(1-q)^{-1/2} = 1 + \frac{3}{2}p + \frac{1}{2}q + R_q$$

и уравнения наши приведутся к виду:

$$\left. \begin{array}{l} \dfrac{d^2p}{d\tau^2} + 2(1+m)i\dfrac{dp}{d\tau} - \dfrac{3}{2}l(1+p) = \dfrac{3}{2}\lambda(q-1)e^{-2i\tau} + lR_p, \\ \dfrac{d^2q}{d\tau^2} - 2(1+m)i\dfrac{dq}{d\tau} - \dfrac{3}{2}l(1+p) = \dfrac{3}{2}\lambda(p-1)e^{+2i\tau} + lR_q. \end{array} \right\} \quad (15)$$

Для этих уравнений мы теперь будем искать периодическое решение с периодом $2\pi$[7], подчиненное требованию, чтобы функции $p$ и $q$ при бесконечно малом $\lambda$ делались также бесконечно малыми, и прежде всего покажем, что всякий раз, когда модуль параметра $\lambda$ достаточно мал, функции $p$ и $q$ в этом решении можно определять рядами, расположенными по целым положительным степеням $\lambda$, и что ряды эти могут содержать только одно произвольное постоянное.

Воспользуемся для этого методом Пуанкаре.

Пусть $p_0, q_0, p'_0, q'_0$ суть начальные значения функций

$$p, \quad q, \quad p' = \frac{dp}{d\tau}, \quad q' = \frac{dq}{d\tau},$$

которые, допустим, соответствуют $\tau = 0$.

Всякий раз, когда модули величин

$$p_0, q_0, p'_0, q'_0, \lambda \quad (16)$$

достаточно малы, функции $p$ и $q$, удовлетворяющие уравнениям (15), и их производные $p'$ и $q'$, на основании одного общего предложения, можно представлять рядами, расположенными по целым положительным степеням величин (16), для всех значений $\tau$, лежащих между 0 и некоторым пределом $T$, который уменьшением модулей величин (16) может быть сделан сколь угодно большим.

Допустим, что таким образом мы сделали $T > 2\pi$. Тогда названными сейчас рядами можем воспользоваться для составлении значений функций $p, q, p', q'$, соответствующих $\tau = 2\pi$.

---

[7] Как увидим, для решения этого не только $2\pi$ но и $\pi$ будет периодом.

Пусть

$$\bar{p}, \bar{q}, \bar{p}', \bar{q}'$$

суть эти значения.

Составляя при посредстве их уравнения

$$\bar{p} = p_0, \bar{q} = q_0, \bar{p}' = p_0', \bar{q}' = q_0', \qquad (17)$$

мы получим условия, необходимые и, очевидно, достаточные для того, чтобы функции $p$ и $q$, определяемые уравнениями (15), были периодическими периода $2\pi$.

Вопрос приводится таким образом к исследованию уравнений (17), которым мы должны удовлетворить выбором постоянных $p_0, q_0, p_0', q_0'$.

Нетрудно показать, что уравнения (17) не все различны и что одно из них есть следствие трех остальных.

Для этого замечаем, что уравнения (15) допускают следующее интегральное уравнение:

$$[p' - i(1-p)][q' + i(1-q)] - \frac{3}{2}m^2(1-p)(1-q) - \frac{2l}{\sqrt{(1-p)(1-q)}}$$
$$= \frac{3}{4}\lambda\big[(1-p)^2 e^{2i\tau} + (1-q)^2 e^{-2i\tau}\big] - \frac{C}{a^3},$$

представляющее собой преобразование уравнения, соответствующего (4) для уравнений (12).

Разлагая входящую сюда функцию

$$(1-p)^{-1/2}(1-q)^{-1/2}$$

в ряд по степеням $p$ и $q$, мы приведем это интегральное уравнение к виду

$$2(1+m)(p+q) - i(p' - q') + \cdots + \frac{3}{4}\lambda\big[(1-p)^2 e^{2i\tau} + (1-q)^2 e^{-2i\tau}\big] = \text{const},$$

где не выписанные члены—выше второго измерения относительно $p, q, p', q'$, и отсюда, подставляя вместо $p, q, p', q'$ их выражении под видом рядов, затем делая $\tau = 2\pi$, получим следующее тождество:

$$2(1+m)(\bar{p} + \bar{q}) - i(\bar{p}' - \bar{q}') + \cdots + \frac{3}{4}\lambda\big[(1-\bar{p})^2 + (1-\bar{q})^2\big]$$
$$= 2(1+m)(p_0 + q_0) - i(p_0' - q_0') + \cdots + \frac{3}{4}\lambda\big[(1-p_0)^2 + (1-q_0)^2\big]$$

относительно величин (16).

Отсюда ясно, что если из уравнений (17) мы удовлетворим каким-либо трем, то четвертое удовлетворится само собой.

На этом основании мы можем ограничиться рассмотрением трех следующих уравнений:

$$\bar{p} = p_0, \quad \bar{p}' = p_0', \quad \bar{q}' = q_0'. \qquad (18)$$

Посмотрим, что могут дать эти уравнения

Обращаясь к уравнениям (15), мы замечаем, что если в них отбросить все члены выше первого измерения относительно $p$ и $q$, а также члены, зависящие от $\lambda$, то уравнения эти обратятся в линейные с постоянными коэффициентами, для которых общий интеграл представится под видом

$$p = \frac{3}{2}lC_1 e^{-\chi i\tau} - \left[\chi^2 + 2(1+m)\chi + \frac{3}{2}l\right]C_2 e^{-\chi i\tau} + C_3\tau + C'$$
$$q = \frac{3}{2}lC_2 e^{-\chi i\tau} - \left[\chi^2 + 2(1+m)\chi + \frac{3}{2}l\right]C_1 e^{+\chi i\tau} - C_3\tau + C'',$$

где

$$\chi = \sqrt{1 + 2m - \frac{1}{2}m^2},$$

а $C_1, C_2, C_3, C', C''$ произвольные постоянные связанные соотношением

$$C' + C'' = \frac{4}{3}\frac{m+1}{l}iC_3.$$

Выражая эти постоянные через начальные значения функций $p, q, p', q'$, мы представим их под видом некоторых линейных комбинаций величин $p_0, q_0, p'_0, q'_0$, причем $C_1$ $C_2$ и $C_3$ выразятся только при помощи трех следующих постоянных:

$$p_0 + q_0, \quad p'_0, \quad q'_0, \tag{19}$$

относительно которых будут независимыми линейными формами.

Разумея теперь под

$$C_1, C_2, C_3 \quad \text{и} \quad C' - C'' = C$$

определяемые таким путем линейные формы величин

$$p_0 + q_0, \quad p'_0, \quad q'_0 \quad \text{и} \quad p_0 - q_0,$$

мы введем эти формы в качестве произвольных постоянных и в общий интеграл точных уравнений (15), который представится таким образом рядами, расположенными по целым положительным степеням величин

$$C_1, C_2, C_3, C \quad \text{и} \quad \lambda.$$

Мы получим при этом для $p$ и $q$ выражения, которые будут отличаться от только что рассмотренных дишь членами выше первого измерения относительно $C_1, C_2, C_3, C$ и $\lambda$ членами, зависящими от $\lambda$.

Вследствие этого, составляя уравнения (18) и выписывая в них только члены, не зависящие от $\lambda$ и притом линейные относительно $C_1, C_2, C_3, C$, найдем, что уравнения эти будут вида

$$\frac{3}{2}l(\rho - 1)C_1 - \left[\chi^2 + 2(1+m)\chi + \frac{3}{2}l\right]\left(\frac{1}{\rho} - 1\right)C_2 + 2\pi C_3 + \cdots = 0,$$
$$\frac{3}{2}l(\rho - 1)C_1 + \left[\chi^2 + 2(1+m)\chi + \frac{3}{2}l\right]\left(\frac{1}{\rho} - 1\right)C_2 + \cdots = 0,$$
$$\left[\chi^2 + 2(1+m)\chi + \frac{3}{2}l\right](\rho - 1)C_1 + \frac{3}{2}l\left(\frac{1}{\rho} - 1\right)C_2 + \cdots = 0,$$

где

$$\rho = e^{2\pi\chi i}$$

и так как функциональный определитель первых их частей в отношении величин $C_1, C_2, C_3$ после положения

$$C_1 = C_2 = C_3 = C = \lambda = 0$$

обращается в величину

$$2\pi(\rho - 1)\left(\frac{1}{\rho} - 1\right)\left\{\frac{3}{4}l^2 - \left[\chi^2 + 2(1+m)\chi + \frac{3}{2}l\right]^2\right\},$$

или, что все равно,

$$-16\pi(1+m)\chi[\chi + 2(1+m)]^2\sin^2\pi\chi,$$

которая при значении $m$, соответствующем лунной теории, не есть нуль, то на основании известного предложения уравнения наши будут разрешимы относительно $C_1, C_2, C_3$ и при условии, что эти последние должны уничтожаться при $C = \lambda = 0$, дадут для них вполне определенные выражении под видом рядов, расположенных по целым положительным степеням $C$ и $\lambda$, абсолютно сходящихся, пока модули $C$ и $\lambda$ достаточно малы.

Так как в получаемых таким путем выражениях $C_1, C_2$, и $C_3$ постоянное $C$ будет фигурировать только в членах выше первого измерения относительно $C$ и $\lambda$, то выражения эти позволят определить величины (19) под видом рядов, расположенных по целым положительным степеням величин $p_0 - q_0$ и $\lambda$.

Таким образом, начальные значения $p_0, q_0, p'_0, q'_0$, а следовательно и самые функции $p, q, p', q'$, соответствующие искомому периодическому решению, представятся рядами, расположенными по целым положительным степеням величин $p_0 - q_0$ и $\lambda$ и будут содержать одно произвольное постоянное $p_0 - q_0$.

Мы будем далее рассматривать это периодическое решение в предположении

$$p_0 = q_0, \qquad (20)$$

равносильном допущению, что при $\tau = 0$ $y = 0$.

Так как $\tau$ уже содержит произвольное постоянное $t_0$, то допущение это не внесет никакого существенного ограничения нашей задачи и послужит только к определению того, что мы будем разуметь под $t_0$, которое представит таким образом один из моментов соединений.

4. На основании изложенного выше мы будем искать наше периодическое решение под видом рядов

$$\left.\begin{aligned} p &= p_1\lambda + p_2\lambda^2 + p_3\lambda^3 + \cdots, \\ q &= q_1\lambda + q_2\lambda^2 + q_3\lambda^3 + \cdots, \end{aligned}\right\} \quad (21)$$

расположенных по целым положительным степеням $\lambda$.

Здесь $p_j, q_j$, суть не зависящие от $\lambda$ периодические функции $\tau$, которые, согласно условию (20), мы будем определять в предположении, что

$$\text{для} \quad \tau = 0 \quad p_j = q_j. \quad (22)$$

Подставляя ряды (21) в уравнения (15) и выражая, что последние удовлетворяются независимо от $\lambda$, получим прежде всего следующую систему уравнений:

$$\left.\begin{aligned} \frac{d^2 p_1}{d\tau^2} + 2(1+m)i\frac{dp_1}{d\tau} - \frac{3}{2}l(p_1 + q_1) &= -\frac{3}{2}e^{-2i\tau}, \\ \frac{d^2 q_1}{d\tau^2} - 2(1+m)i\frac{dq_1}{d\tau} - \frac{3}{2}l(p_1 + q_1) &= -\frac{3}{2}e^{+2i\tau}, \end{aligned}\right\} \quad (23)$$

которой определятся $p_1$ и $q_1$. Затем, если допустим, что разложение по восходящим степеням $\lambda$ после подстановки выражений (21) дает

$$\left.\begin{aligned} R_p &= P_2\lambda^2 + P_3\lambda^3 + P_4\lambda^4 \ldots, \\ R_q &= Q_2\lambda^2 + Q_3\lambda^3 + Q_4\lambda^4 \ldots, \end{aligned}\right\}$$

получим ряд систем уравнений вида

$$\left.\begin{aligned} \frac{d^2 p_j}{d\tau^2} + 2(1+m)i\frac{dp_j}{d\tau} - \frac{3}{2}l(p_j + q_j) &= \frac{3}{2}q_{j-1}e^{-2i\tau} + lP_j, \\ \frac{d^2 q_j}{d\tau^2} - 2(1+m)i\frac{dq_j}{d\tau} - \frac{3}{2}l(p_j + q_j) &= \frac{3}{2}p_{j-1}e^{+2i\tau} + lQ_j, \end{aligned}\right\} \quad (24)$$

из которых при $j = 2,3 \ldots$ последовательно найдем $p_2$ и $q_2$, $p_3$ и $q_3$ и т. д., так как $P_j$ и $Q_j$, очевидно, будут зависеть только от величин

$$p_1, p_2, \ldots, p_{j-1}, q_1, q_2, \ldots, q_{j-1}. \quad (25)$$

Этими уравнениями при условии периодичности и при условии (22) функции $p_j, q_j$ определятся вполне и, как нетрудно убедиться, будут вида

$$\left.\begin{aligned} p_j &= \sum_{s=0}^{s=j} a_{j,j-2s} e^{2(j-2s)i\tau}, \\ q_j &= \sum_{s=0}^{s=j} a_{j,j-2s} e^{-2(j-2s)i\tau}, \end{aligned}\right\} \quad (26)$$

где $a_{j,\sigma}$, суть некоторые постоянные.

Действительно, из уравнений (23) при сказанных условиях находим

$$p_1 = a_{1,1}e^{2i\tau} + a_{1,-1}e^{-2i\tau}, \quad q_1 = a_{1,1}e^{-2i\tau} + a_{1,-1}e^{2i\tau},$$

где

$$a_{1,1} = -\frac{9}{16}\frac{2 + 4m + 3m^2}{6 - 4m + m^2}, \quad a_{1,-1} = \frac{3}{16}\frac{38 + 28m + 9m^2}{6 - 4m + m^2},$$

а допуская, что выражения (26) справедливы для всех значений $j$, меньших какого-либо числа $r$, легко докажем, что они будут справедливы и для $j = r$.

Для этого замечаем, что $P_j$, и $Q_j$, как это следует из самого их определения, суть полиномы, составленные из величин (25) так, что вес каждого члена равен $j$. Поэтому, если, например, в выражении

$$\frac{3}{2}q_{r-1}e^{-2i\tau} + lP_r$$

все входящие в него $p_j$ и $q_j$ заменим выражениями (26), которые можно представить так:

$$p_j = e^{2ji\tau}\sum_{s=0}^{s=j}a_{j,j-2s}e^{-4si\tau}, \quad q_j = e^{2ji\tau}\sum_{s=0}^{s=j}a_{j,2s-j}e^{-4si\tau},$$

то получим результат вида

$$\frac{3}{2}q_{r-1}e^{-2i\tau} + lP_r = \sum_{s=0}^{s=j}A_{r,r-2s}e^{2(r-2s)i\tau},$$

где $A_{r,r-2s}$ суть постоянные, известным образом зависящие от постоянных $a_{j,\sigma}$, соответствующих $j < r$.

Вместе с тем будем иметь

$$\frac{3}{2}p_{r-1}e^{2i\tau} + lQ_r = \sum_{s=0}^{s=j}A_{r,r-2s}e^{-2(r-2s)i\tau},$$

ибо $Q_r$ получается из $P_r$, перестановкой букв $p$ и $q$, а в формулах (26), $q_j$ выводится из $p_j$ заменой $\tau$ на $-\tau$.

Вследствие этого, если все $p_j$, и $q_j$, для которых $j < r$, представляются выражениями (26), то уравнения (24) для $j = r$ примут вид

$$\left.\begin{aligned}\frac{d^2p_r}{d\tau^2} + 2(1+m)i\frac{dp_r}{d\tau} - \frac{3}{2}l(p_r + q_r) &= \sum_{s=0}^{s=j}A_{r,r-2s}e^{2(r-2s)i\tau},\\ \frac{d^2q_r}{d\tau^2} - 2(1+m)i\frac{dq_r}{d\tau} - \frac{3}{2}l(p_r + q_r) &= \sum_{s=0}^{s=j}A_{r,r-2s}e^{-2(r-2s)i\tau},\end{aligned}\right\}$$

и при поставленных нами требованиях дадут для $p_r$ и $q_r$ выражения того же типа (26), причем коэффициенты $a_{r,\sigma}$ определятся уравнениями

вида

$$\left[4\sigma^2 + 4(1+m)\sigma + \frac{3}{2}l\right]a_{r,\sigma} + \frac{3}{2}la_{r,-\sigma} = -A_{r,\sigma} \ (\sigma = -r, -r+2, \ldots, r-2, r),$$

из которых найдем

$$a_{r,\sigma} = \frac{\frac{3}{2}lA_{r,-\sigma} - \left[4\sigma^2 + 4(1+m)\sigma + \frac{3}{2}l\right]A_{r,\sigma}}{2\sigma^2[2(4\sigma^2 - 1) - 4m + m^2]}, \qquad (27)$$

если $\sigma$ не нуль, и

$$a_{r,0} = -\frac{1}{3l}A_{r,0}. \qquad (28)$$

Последняя формула относится к случаю четного $r$.

Таким образом, зная, что $p_1$ и $q_1$ определяются выражениями вида (26), мы можем быть уверены, того же вида выражения получатся и для всех остальных $p_j$, и $q_j$. При этом, зная $a_{1,1}$ и $a_{1,-1}$ можем вычислять все остальные коэффициенты по формулам (27) и (28), которые дают все $a_{r,\sigma}$, соответствующе какому-либо $r$, когда известны все $a_{j,s}$, для которых $j < r$.

Таким образом убеждаемся, что искомые функции $p$ и $q$ представятся и рядами

$$\left.\begin{array}{l} p = \sum\limits_{j=1}^{j=\infty} \sum\limits_{s=0}^{s=j} a_{j,j-2s}\lambda^j e^{2(j-2s)i\tau}, \\ q = \sum\limits_{j=1}^{j=\infty} \sum\limits_{s=0}^{s=j} a_{j,j-2s}\lambda^j e^{2(2s-j)i\tau}, \end{array}\right\} \qquad (29)$$

откуда для функций $\xi$ и $\eta$ получим следующие выражения:

$$\left.\begin{array}{l} \xi = -\sum\limits_{j=1}^{j=\infty} \sum\limits_{s=0}^{s=j} a_{j,j-2s}\lambda^j \cos 2(j-2s)\tau, \\ \eta = -\sum\limits_{j=1}^{j=\infty} \sum\limits_{s=0}^{s=j} a_{j,j-2s}\lambda^j \sin 2(j-2s)\tau, \end{array}\right\} \qquad (30)$$

Располагая ряды (29) по степеням $e^{2i\tau}$ и делая

$$a(1-p) = \sum a_s e^{2si\tau}, \quad a(1-q) = \sum a_s e^{-2si\tau},$$

будем иметь

$$\left.\begin{array}{l} a_0 = a\left\{1 - \sum\limits_{s=1}^{s=\infty} a_{2s,0}\lambda^{2s}\right\}, \\ a_{\pm\sigma} = -a\sum\limits_{s=0}^{s=\infty} a_{\sigma+2s,\pm\sigma}\lambda^{\sigma+2s}, \end{array}\right\} \qquad (31)$$

где под $\sigma$ разумеется число положительное.

Отсюда для функций $u$ и $v$ получаем выражения (2), в которых, согласно допущению Хилла, коэффициенты $a_\sigma$ таковы, что при $\lambda$ бесконечном малом первого порядка отношение $\frac{a_{\pm\sigma}}{a_0}$ есть бесконечно малое порядка $\sigma$.

Формулы (31) дают выражения коэффициентов $a_s$ рассматриваемых Хиллом под видом рядов, расположенных по степеням $\lambda$. Чтобы получить такие же выражения из формул Хилла, мы должны были бы формулу (11) заменить следующей:

$$a_0 = aJ^{1/3}\left\{1 + \frac{3}{2}\frac{\lambda}{l}\frac{a_{-1}}{a_0}\right\}^{-1/3} \tag{32}$$

и это выражение $a_0$, а затем выражения

$$a_{\pm\sigma} = a_0\lambda^\sigma \sum_{s=0}^{s=\infty} \alpha_{\pm\sigma,s}\lambda^{2\sigma} \tag{33}$$

разложить в ряды по степеням $\lambda$.

Мы получили бы таким путем формулы, дающие выражения коэффициентов $a_{r,\sigma}$ через коэффициенты $\alpha_{j,s}$, вычисляемые по формулам Хилла, как было показано в п. 2, и должно заметить, что для вычисления формулы эти удобнее наших. Но из них весьма трудно выводить какие-либо заключения о сходимости рассматриваемых здесь рядов, тогда как наши формулы (27) и (28) весьма удобны для этой цели.

Обращаясь теперь к вопросу о сходимости, мы будем рассматривать ряды (29) как расположенные по степеням $\lambda$ и $e^{2i\tau}$ и при этом, ограничиваясь вещественными значениями $\tau$, покажем, что для величин $m$ и $\lambda$, соответствующих лунной теории, ряды эти суть абсолютно сходящиеся, откуда будет следовать и сходимость при тех же условиях рядов (30) и (31).

Мы увидим, что при сказанных условиях ряды, составленные из модулей членов рядов (29), будут менее 1. Поэтому и подавно будет иметь место неравенство

$$\sum_{s=1}^{s=\infty}|a_{2s,0}|\lambda^{2s} < 1,$$

если вообще под $|a|$ будем разуметь модуль количества $a$.

Отсюда, вследствие абсолютной сходимости рядов (31), заключим о такой же сходимости рядов

$$\alpha_{\pm\sigma,0}\lambda^\sigma + \alpha_{\pm\sigma,1}\lambda^{\sigma+2} + \alpha_{\pm\sigma,2}\lambda^{\sigma+4} + \cdots,$$

которыми мы представили отношения $\frac{a_{\pm\sigma}}{a_0}$.

Мы рассмотрим также ряды, в которые обращаются ряды (29) после замены коэффициентов $a_{j,r}$ их разложениями по степеням $m$.

Эти коэффициенты, как показывают формулы (27) и (28), будут рациональными функциями $m$ и представятся дробями, знаменатели которых будут содержать только множители типа

$$2(4\sigma^2 - 1) - 4m + m^2, \qquad (\sigma = 1, 2, 3, \dots)$$

и

$$l = 1 + 2m + \frac{3}{2}m^2.$$

Поэтому при $|m| < \sqrt{\frac{2}{3}}$ они будут разложимы в ряды по целым положительным степеням $m$.

Предполагая, что эти разложения подставлены в ряды (29), и рассматривая полученные из них таким путем новые ряды, как расположенные по степеням $m, \lambda$ и $e^{2i\tau}$, мы покажем, что при величинах $m$ и $\lambda$, соответствующих лунной теории, и при вещественных значениях $\tau$ и эти новые ряды будут абсолютно сходящимися. Поэтому и подавно будут абсолютно сходящимися ряды, расположенные по степеням $m, e^{2i\tau}$, в которые обратятся сейчас указанные, если в них сделаем $\lambda = m^2$.

Мы докажем таким образом сходимость рядов, расположенных по степеням $m$, которыми могут быть представлены коэффициенты $a_s$.

Наконец, если, вместо $m$, пожелаем ввести величину $m_1 = \frac{m}{1+m}$, представляющую отношение среднего движения Солнца к среднему звездному движению Луны, и, подобно тому как это делалось в некоторых старых теориях Луны, коэффициенты $a_s$, представим рядами, расположенными по степеням $m_1$, то сходимость этих рядов будет следовать тотчас же из абсолютной сходимости рядов, расположенных по степеням $m$, если принять в расчет, что разложение $m$ по степеням $m_1$, обладает положительными коэффициентами.

5. Предполагая $m$ числом положительным и рассматривая формулу (27) при различных значениях $\sigma$ взятых из ряда

$$\pm 1, \pm 2, \pm 3, \dots,$$

замечаем тотчас же, что коэффициенты при $A_{r,\sigma}$ и $A_{r,-\sigma}$ достигают наибольших численных значений в случае $\sigma = -1$. Мы будем поэтому иметь неравенство

$$|a_{r,\sigma}| < \frac{\frac{3}{2}l|A_{r,-\sigma}| + \left(8 + 4m + \frac{3}{2}l\right)|A_{r,\sigma}|}{2(6 - 4m + m^2)}, \qquad (34)$$

справедливое для всех указанных сейчас значений $\sigma$.

Но неравенство — это справедливо и для $\sigma = 0$.

Действительно, в этом случае вторая часть его приводится к

$$\frac{8 + 4m + 3l}{2(6 - 4m + m^2)}|A_{r,0}|,$$

что представляет величину, большую

$$|a_{r,0}| = \frac{1}{3l}|A_{r,0}|,$$

ибо при всяком положительном $m$ разность

$$\frac{8+4m+3l}{2(6-4m+m^2)} - \frac{1}{3l} = \frac{22+20m+9m^2}{4(6-4m+m^2)} - \frac{2}{6+12m+9m^2},$$

очевидно, положительна.

Таким образом, неравенство (34) справедливо для всех значений $\sigma$, которые здесь приходится рассматривать.

Исходя из этого неравенства, нетрудно теперь получить формулу, дающую высшие пределы модулей величин $a_{r,\sigma}$, соответствующих какому-либо $r$, когда известны высшие пределы модулей всех $a_{j,\sigma}$, для которых $j < r$.

Для этого замечаем, что $A_{r,\sigma}$, есть полином, составленный из величин $a_{j,s}$, в котором все коэффициенты положительны. Это вытекает из самого определения $A_{r,\sigma}$, если принять в расчет, что коэффициенты разложений $R_p$ и $R_q$ по степеням $p$ и $q$ все положительны.

Отсюда следует, что, заменяя в выражении $A_{r,\sigma}$ все $a_{j,s}$ их модулями или какими-либо высшими пределами их модулей, мы получим некоторый высший предел для $|A_{r,\sigma}|$.

Допустим, что каким-либо способом мы нашли высшие пределы модулей для всех $a_{j,s}$, для которых $j < r$. Означая эти высшие пределы через $\mathbf{a}_{j,s}$ и результат замены в выражении $A_{r,\sigma}$ величин $a_{j,s}$ величинами $\mathbf{a}_{j,s}$ через $\mathbf{A}_{r,\sigma}$ положим

$$\mathbf{a}_{r,\sigma} = \frac{\frac{3}{2}l\mathbf{A}_{r,-\sigma} + \left(8+4m+\frac{3}{2}l\right)\mathbf{A}_{r,\sigma}}{2(6-4m+m^2)^2}. \tag{35}$$

Тогда на основании замеченного сейчас это выражение $\mathbf{a}_{r,\sigma}$ будет некоторым высшим пределом модуля $a_{r,\sigma}$.

Таким образом, принимая

$$\mathbf{a}_{1,1} = |a_{1,1}|, \qquad \mathbf{a}_{1,-1} = |a_{1,-1}|$$

и определяя все остальные $\mathbf{a}_{j,s}$ по формуле (36), мы цолучим высшие пределы модулей всех $a_{j,s}$.

Нетрудно притом видеть, что, если принять

$$\mathbf{A}_{1,-1} = \frac{3}{2}, \qquad \mathbf{A}_{1,1} = 0,$$

формула (35) при $r = 1$ будет приводить к равенствам (36) и, следовательно, ею можно будет пользоваться для всех значений $r$, которые необходимо здесь рассматривать,

Заметим, что определяемые таким путем величины $\mathbf{a}_{j,s}$ все будут возрастающими функциями $m$, пока $m$, оставаясь положительным, не превос-ходит некоторого предела. Так, например, это, очевидно, будет иметь место, пока $m < 2$ ибо при этом выражение

$$6 - 4m + m^2$$

будет убывающей функцией $m$

Рассмотрим теперь ряд

$$\sum_{r=1}^{r=\infty}\sum_{s=0}^{s=r}\mathbf{a}_{r,r-2s}\lambda^r, \tag{37}$$

и допустим, что для каких-либо положительных значений $m$ и $\lambda$ которые назовем $m_1$ и $\lambda_1$, доказана его сходимость. Тогда, если $m_1 < 2$, будет доказана абсолютная сходимость рядов (29) при вещественном $\tau$ и при $|\lambda| < \lambda_1$ для всякого положительного $m$, которое не превосходит $m_1$.

Но условии сходимости ряда (37) найти весьма нетрудно, ибо ряд этот можно определить некоторым алгебраическим уравнением.

Действительно, из (35) следует

$$\sum_{s=0}^{s=r}\mathbf{a}_{r,r-2s} = \frac{8 + 4m + 3l}{2(6 - 4m + m^2)}\sum_{s=0}^{s=r}\mathbf{A}_{r,r-2s},$$

а сумма, находящаяся во второй части равенства, есть результат замены в выражении

$$\frac{3}{2}q_{r-1} + lP_r$$

каждого $p_j$, и каждого $q_j$ выражением

$$\sum_{s=0}^{s=j}\mathbf{a}_{j,j-2s}$$

соответствующим тому, же $j$.

Поэтому если ряд (37) рассматривать как простой, общий член которого представляется выражением

$$\sum_{s=0}^{s=r}\mathbf{a}_{r,r-2s}\lambda^r,$$

то ряд этот будет получаться при разложении по степеням $\lambda$ уничтожающеюся при $\lambda = 0$ корня $z$ уравнения

$$z = \frac{8 + 4m + 3l}{2(6 - 4m + m^2)}\left\{\frac{3}{2}(1 + z)\lambda + l\left[\frac{1}{(1-z)^2} - 1 - 2z\right]\right\}.$$

Принимая в расчет выражение $l$ и полагая

$$\varepsilon = \frac{3(22 + 20m + 9m^2)}{8(6 - 4m + m^2)}\lambda, \quad h = \frac{22 + 20m + 9m^2}{4(6 - 4m + m^2)}l, \tag{38}$$

уравнение это приведем к виду

$$z = \varepsilon(1+z) + h\frac{(3-2z)z^2}{(1-z)^2}, \tag{39}$$

и задача наша приведется таким образом к разысканию условия сходимости ряда, расположенного по степеням $\varepsilon$, в который, при достаточно малом $\varepsilon$, разлагается корень уравнения (39), уничтожающийся при $\varepsilon = 0$.

Мы рассматривали до сих пор только положительные значения $m$. Будем теперь рассматривать всякие его значения, как вещественные, так и комплексные, модули которых не превосходят некоторого предела $M$, и покажем, что при достаточно малом $M$ формула (35), если в ней $m$ заменить через $M$, будет давать высшие пределы модулей величин $a_{r,\sigma}$, при всяком $m$, модуль которого не превосходит $M$.

Мы будем основываться при этом на следующем легко доказываемом предложении:

Если $a, b, c$ суть положительные числа, удовлетворяющие неравенству

$$ac - b^2 > 0,$$

а $M$ — положительное число, не превосходящее меньшего корня уравнения

$$\frac{x^2}{a} - 2\frac{x}{b} + \frac{1}{c} = 0,$$

то низшим пределом модуля функции

$$a - 2bx + cx^2$$

при величинах $x$, модули которых не превосходят $M$, будет число

$$a - 2bM + cM^2.$$

Применяя это предложение к функции

$$2(4\sigma^2 - 1) - 4m + m^2,$$

входящей в знаменатель выражения (27), приходим к заключению, что, если $M$ не превосходит числа

$$\frac{2\sqrt{4\sigma^2 - 1}}{\sqrt{4\sigma^2 - 1} + \sqrt{4\sigma^2 - 3}} \tag{40}$$

(где радикалы предполагаются положительными), то при $|m| \leq M$ будет выполняться неравенство

$$|2(4\sigma^2 - 1) - 4m + m^2| \geq 2(4\sigma^2 - 1) - 4M + M^2.$$

Так как число (40), всегда большее единицы, с беспредельным возрастанием $\sigma$ стремится к пределу, равному 1, то, желая, чтобы написанное сейчас неравенство выполнялось для всех значений $\sigma$, которые приходится рассматривать, имея дело с формулой (27), мы должны предположить $M \leq 1$. В этом предположении будем иметь

$$|2(4\sigma^2 - 1) - 4m + m^2| \geq 6 - 4M + M^2$$

для всех значений $m$, модули которых не превосходят $M$.

Отсюда заключаем, что если при сделанном предположении в формуле (35) $m$ заменим через $M$ (как в коэффициентах при $\mathbf{A}_{r,-\sigma}$, и $\mathbf{A}_{r,\sigma}$, так и в выражениях последних, которые содержат $l$) и под обозначениями $\mathbf{a}_{j,s}$, входящими в выражения $\mathbf{A}_{r,-\sigma}$, и $\mathbf{A}_{r,\sigma}$, будем разуметь какие-либо высшие пределы модулей величин $a_{j,s}$, при $|m| \leq M$, то по крайней мере для $\sigma$ отличного от нуля формула эта будет давать высшие пределы модулей величин $a_{r,s}$, при том же условии $|m| \leq M$.

Но нетрудно что при достаточно малом $M$ формула эта будет давать также и высший предел модуля $a_{r,0}$.

Действительно, на основании указанного выше предложения модуль функции

$$l = 1 + 2m + \frac{3}{2}m^2$$

при $|m| \leq M$ будет не менее

$$1 - 2M + \frac{3}{2}M^2$$

всякий раз, когда $M$ не превосходит числа

$$1 - \frac{1}{\sqrt{3}} = 0{,}42 \ldots$$

Поэтому, предполагая

$$M \leq 1 - \frac{1}{\sqrt{3}}, \tag{41}$$

из (28) найдём

$$|a_{r,0}| < \frac{2\mathbf{A}_{r,0}}{3(2 - 4M + 3M^2)},$$

вследствие чего, замечая, что разность

$$\frac{22 + 20M + 9M^2}{4(6 - 4M + M^2)} - \frac{2}{3(2 - 4M + 3M^2)} = \frac{84 - 112M + 4M^2 + 72M^3 + 81M^4}{12(6 - 4M + M^2)(2 - 4M + 3M^2)}$$

положительна при всяком вещественном $M$, и определяя по формуле (35) (после вышеуказанной замены) $\mathbf{a}_{r,0}$, заключим, что

$$|a_{r,0}| < \mathbf{a}_{r,0}.$$

Таким образом, убеждаемся, что если при условии (41) в формуле (35) $m$ заменим через $M$ и формулой этой будем пользоваться, начиная от $r = 1$, принимая по-прежнему

$$\mathbf{A}_{1,-1} = \frac{3}{2}, \qquad \mathbf{A}_{1,1} = 0,$$

то она будет давать высшие пределы для модулей величин $a_{r,\sigma}$, каково бы ни было $m$, модуль которого не превосходит $M$.

Допустим теперь, что рассматриваются ряды, расположенные по степеням $m, \lambda,$ и $e^{2i\tau}$, в которые обращаются ряды (29) после разложения коэффициентов $a_{r,\sigma}$, по восходящим степеням $m$.

Принимая условие (41) и разумея под $\mathbf{a}_{r,\sigma}$, высшие пределы, о которых сейчас шла речь, на основании известной теоремы заключим, что коэффициенты разложения $a_{r,\sigma}$, по восходящим степеням $m$ будут численно менее соответственных коэффициентов такого же разложения функции

$$\mathbf{a}_{r,\sigma}\left(1-\frac{m}{M}\right)^{-1}.$$

Поэтому для доказательства сходимости рассматриваемых рядов при вещественном $\tau$ и при условиях

$$|m| < M, \qquad |\lambda| \leq \lambda_1$$

достаточно будет доказать сходимость ряда (37) при $\lambda = \lambda_1$.

Таким образом, вопрос о сходимости этих рядов приведется к решению той самой задачи, к которой мы пришли, рассматривая ряды (29) при $m$ положительном.

6. Обращаемся к уравнению (39), которое будем рассматривать в предположении, что $\varepsilon$ и $h$ суть положительные числа.

Задача наша привилась к определению условии сходимости ряда

$$\varepsilon + c_2\varepsilon^2 + c_3\varepsilon^3 + \cdots, \tag{42}$$

расположенного до восходящим степеням $\varepsilon$, в который при достаточно малом $\varepsilon$ разлагается корень уравнения (39), уничтожающийся при $\varepsilon = 0$.

Укажем сначала одно весьма простое достаточное для этой сходимости условие.

Это условие получим тотчас же, если, вместо уравнения (39), которое третьей степени относительно $z$, будем рассматривать следующее квадратное уравнение:

$$z = \varepsilon(1+z) + \frac{9hz^2}{3-4z}. \tag{43}$$

Нетрудно заметить, что коэффициенты разложения по восходящим степеням $z$ второй части этого последнего уравнения не менее соответственных коэффициентов такого же разложения второй части уравнения (39). Поэтому ряд вида (42), в который при достаточно малом $\varepsilon$ разложится корень уравнения (43), уничтожающийся при $\varepsilon = 0$, будет обладать не меньшими коэффициентами, чем ряд (42), и для сходимости последнего будет достаточно, чтобы был сходящимся этот новый ряд.

Но названный сейчас ряд получается при разложении по восходящим степеням $\varepsilon$ функции

$$z = \frac{3 + \varepsilon - \sqrt{9 - 6(7 + 18h)\varepsilon + 49\varepsilon^2}}{2(4 + 9h - 4\varepsilon)},$$

рассмотрение которой приводит к следующему условию его сходимости:

$$\varepsilon \leq \frac{3}{7 + 18h + \sqrt{(7 + 18h)^2 - 49}},$$

где радикал должно считать положительным.

Отсюда видно, что, если

$$\varepsilon \leq \frac{3}{2(7 + 18h)},$$

а тем более, если

$$\varepsilon \leq \frac{1}{6(1 + 2h)}, \tag{44}$$

рассматриваемый ряд, а следовательно, и ряд (42), будет сходящимся.

Условие (44) мы и желали указать.

Не останавливаясь на доказательстве, скажем теперь, в чем состоит условие, не только достаточное для сходимости ряда (42), но и необходимое.

Положим

$$\omega = \left(\frac{2h}{1 - \varepsilon + 2h}\right)^{1/3}.$$

Тогда это условие выразится так:

$$\varepsilon \leq \frac{2 - \omega - \omega^2}{4 + \omega + \omega^2} < 1, \tag{45}$$

что, очевидно, требует, чтобы при данном $h$ число $\varepsilon$ не превосходило некоторого предела $\varepsilon_0$, лежащего между $0$ и $\frac{1}{2}$.

Можно заметить, что условие это совпадает с тем, при котором уравнение (39) имеет корни между $0$ и $1$, и что ряд (42) представляет наименьший из таких корней.

Обращаемся теперь к выражениям $\varepsilon$ и $h$, которые даются формулами (38), и полагаем

$$\lambda = m^2.$$

При этом $\varepsilon$ и $h$ делаются функциями одного $m$, и возникает вопрос о наибольшем значении $m$, при котором ряд (42) остается сходящимся (будучи сходящимся и для всякого меньшего $m$).

Рассматривая условие (45), можем убедиться, что значение это лежит между $\frac{1}{7}$ и $\frac{1}{6}$. Но уже условия (44) достаточно, чтобы показать, что $\frac{1}{7}$ еще не достигает этого предела.

Действительно, при $m = \frac{1}{7}$ формулы (38) дают

$$\varepsilon = \frac{1227}{34\,888} = \frac{8589}{244\,216}, \quad h = \frac{53\,761}{34\,888},$$

откуда следует

$$\frac{1}{6(1+2h)} = \frac{8\,722}{210\,615} > \varepsilon.$$

Таким образом, убеждаемся, что при $m \leq \frac{1}{7}$ все рассмотренные нами ряды будут сходящимися.

Для теории Луны этого вывода более чем достаточно, так как величина $m$, соответствующая Луне, не превосходит $\frac{1}{12}$.

7. Уравнение (39), послужившее нам для решения вопроса о сходимости рядов (29) и (30), может служить также и для определения высшего предела погрешности, которую делаем, отбрасывая этих рядах все члены выше какого-либо порядка относительно $\lambda$, ибо, отбрасывая в каком-либо из этих рядов все члены выше $n$-го порядка, мы получаем результат, погрешность которого, очевидно, не превосходит величины

$$z - (\varepsilon + c_2\varepsilon^2 + \cdots + c_n\varepsilon^n), \tag{46}$$

где $z$ — наименьший корень уравнения (39).

Нельзя, однако, ожидать, чтобы этот высший предел обладал той степенью точности, какова желательна в рассматриваемой задаче.

Причины его недостаточной точности коренятся, конечно, в самой сущности нашего анализа, и важнейшие из них едва ли легко устранимы. Однако одна из причин неточности может быть легко устранена.

В самом деле, к уравнению (39) мы пришли, рассматривая ряд (37) в предположении, что

$$\mathbf{a}_{1,-1} = |a_{1,-1}|, \quad \mathbf{a}_{1,1} = |a_{1,1}|$$

и что все остальные $\mathbf{a}_{r,\sigma}$, начиная от $r = 2$, определяются формулой (35). Но мы могли бы этой последней формулой пользоваться только при $r > N$, разумея под $N$ какое-либо число, большее 1, и принимать

$$\mathbf{a}_{r,\sigma} = |a_{r,\sigma}|$$

при $r \leq N$. Тогда остаток ряда (37) посыле членов $n$-го порядка относительно $\lambda$ представил бы, очевидно, более точный высший предел погрешности, чем указанное выше число (46).

Рассмотрим же теперь ряд (37) в этом новом предположении.

Нетрудно видеть, что ряд этот по-прежнему можно определять некоторым алгебраическим уравнением третьей степени.

Чтобы составить это уравнение, допустим, что в функцию

$$\varepsilon z + h \frac{(3-2z)z^2}{(1-z)^2}$$

при величинах (38) для $\varepsilon$ и $h$ мы подставили вместо $z$ выражение

$$\sum_{j=1}^{j=N} \sum_{s=0}^{s=j} |a_{j,j-2s}| \lambda^j$$

и, принимая в расчет зависимость $\varepsilon$ от $\lambda$, результат разложили в ряд

$$L_2 \lambda^2 + L_3 \lambda^3 + \cdots + L_N \lambda^N + \cdots$$

по восходящим степеням $\lambda$. Тогда, полагая

$$\varepsilon' = \sum_{j=1}^{j=N} \sum_{s=0}^{s=j} |a_{j,j-2s}| \lambda^j - \sum_{j=2}^{j=N} L_j \lambda^j,$$

искомое уравнение можем представить так:

$$z = \varepsilon' + \varepsilon z + h \frac{(3-2z)z^2}{(1-z)^2}. \tag{47}$$

Допустим теперь, что

$$l_1 \lambda + l_2 \lambda^2 + l_3 \lambda^3 + \cdots \tag{48}$$

есть ряд, в который разлагается по степеням $\lambda$ корень этого уравнения, уничтожающийся при $\lambda = 0$. Тогда погрешность вывода, который получим, останавливаясь, при вычислении какого-либо из рядов (30), на членах $n$-го порядка относительно $\lambda$, не превзойдет величины

$$z - (l_1 \lambda + l_2 \lambda^2 + \cdots + l_n \lambda^n),$$

где $z$ наименьший корень уравнения (47).

Мы получаем таким путем высший предел погрешности, который при $N$ достаточно большом может достигать надлежащей точности. Но возможно, что для этого придется брать вообще $N > n$, как это имеет место в приводимом ниже примере.

Для получения рассматриваемого высшего предела необходимо вычисление наименьшего корня уравнения (47), для чего можно прибегнуть к методе последовательных приближений, приведя предварительно уравнение это к более простому виду

$$z = \delta + g \frac{(3-2z)z^2}{(1-z)^2}, \tag{49}$$

где

$$\delta = \frac{\varepsilon'}{1-\varepsilon}, \quad g = \frac{h}{1-\varepsilon}.$$

В предположении, которое здесь везде подразумевается, что $\lambda, \varepsilon$ и $h$ суть положительные числа, удовлетворяющие условию сходимости ряда (48), можно доказать, что $\delta$ и $g$ будут положительными числами, удовлетворяющими неравенству

$$\frac{1+2\delta}{3(1+2g)} < 1,$$

и что наименьший корень уравнения (49) (которое при сказанном предположении будет иметь только вещественные и положительные корни) не будет превосходить числа

$$\frac{1+2\delta}{3(1+2g)}, \qquad (50)$$

не превосходящего ни одного из двух корней.

При этих условиях, если под $z_0$ будем разуметь какое-либо положительное число, не превосходящее (50), и по этому числу $z_0$ будем определять последовательно $z_1, z_2, z_3$ и т. д., пользуясь уравнением

$$z_{n+1} = \delta + g\frac{(3-2z_n)z_n^2}{(1-z_n)^2},$$

то $z_n$, по мере возрастания $n$, будет стремиться к искомому корню $c$, или постоянно возрастая (когда $z_0 < c$), или постоянно убывая (когда $z_0 > c$).

Таким способом интересующий нас корень можно вычислять довольно быстро с большой точностью.

8. Применим сейчас изложенное к лунной теории, ограничиваясь предположением, что $N = 2$.

В этом предположении, замечая, что

$$L_2\lambda^2 = c_2\varepsilon^2 = (1+3h)\varepsilon^2,$$

и полагая

$$\varepsilon_1 = \{|a_{2,0}| + |a_{2,2}| + |a_{2,-2}|\}\lambda^2,$$

находим

$$\varepsilon' = \varepsilon + \varepsilon_1 - (1+3h)\varepsilon^2,$$
$$l_1\lambda = \varepsilon, \quad l_2\lambda^2 = \varepsilon_1, \quad l_3\lambda^3 = (1+6h)\varepsilon\varepsilon_1 + 4h\varepsilon^3,$$

Для вычисления этих выражений необходимо предварительно найти $a_{2,0}, a_{2,-2}$ и $a_{2,2}$.

Когда известны $a_{1,-1}$ и $a_{1,1}$ которые найдутся по формулам

$$a_{1,-1} = \frac{3}{16}\frac{38+28m+9m^2}{6-4m+m^2}, \quad a_{1,1} = -\frac{9}{16}\frac{2+4m+3m^2}{6-4m+m^2},$$

и вычислено

$$2l = 2 + 4m + 3m^2,$$

$a_{2,0}$ найдется по формуле

$$a_{2,0} = -\frac{1}{4}(a_{1,-1} - a_{1,1})^2 - \left(2a_{1,1} + \frac{1}{2l}\right)a_{1,-1},$$

которая легко выводится как из (28), так и из (31) и (32).

Что касается $a_{2,-2}$ и $a_{2,2}$ то для вычисления их, вместо формулы (27), обращаемся к формулам Хилла, которые скорее приводят к цели.

Из (31) и (33) находим

$$a_{2,-2} = -\alpha_{-2,0}, \qquad a_{2,2} = -\alpha_{2,0},$$

а из формул п. 2

$$\alpha_{2,0} = 2\overline{[2]}\alpha_{1,0} + [2,1]\alpha_{1,0}\alpha_{-1,0},$$
$$\alpha_{-2,0} = 2\overline{(-2)}\alpha_{1,0} + [-2,-1]\alpha_{1,0}\alpha_{-1,0}.$$

Отсюда, принимая в расчет, что

$$\alpha_{-1,0} = -a_{1,-1}, \qquad \alpha_{1,0} = -a_{1,1},$$
$$\overline{[2]} = \frac{3}{64}\frac{2 + 16m + 9m^2}{30 - 4m + m^2}, \qquad \overline{(-2)} = -\frac{9}{64}\frac{38 + 16m + 3m^2}{30 - 4m + m^2},$$
$$[2,1] = -\frac{1}{2}\frac{26 + m^2}{30 - 4m + m^2}, \qquad [-2,-1] = -\frac{1}{2}\frac{18 - 8m + m^2}{30 - 4m + m^2}.$$

выводим

$$a_{2,-2} = \frac{16(18 - 8m + m^2)a_{1,-1} - 9(38 + 16m + 3m^2)a_{1,1}}{32(30 - 4m + m^2)},$$
$$a_{2,2} = \frac{16(26 + m^2)a_{1,-1} - 3(2 + 16m + 9m^2)a_{1,1}}{32(30 - 4m + m^2)}.$$

По этим формулам, принимая $\lambda = m^2$, найдем

$$a_{1,-1}\lambda, \quad a_{1,1}\lambda, \quad a_{2,-2}\lambda^2, \quad a_{2,2}\lambda^2,$$

что тотчас же доставит

$$\varepsilon = \{|a_{1,-1}| + |a_{1,1}|\}\lambda$$

и $\varepsilon_1$, а вычислив затем $h$ по формуле (38), найдем $l_3\lambda^3, \varepsilon', \delta$ и $g$.

Мы примем вместе с Хиллом

$$m = 0{,}08084\,89338\,08312$$

и при выписывании каждого результата вычисления будем держаться обычного требования, чтобы погрешность не превышала половины единицы последнего десятичного знака, на котором останавливаемся. Ограничиваясь при этом десятью десятичными знаками, найдем

$$\begin{aligned}
a_{1,-1}\lambda &= \phantom{-}0{,}00869\,58085\,(-)^8, & h &= 1{,}22011\,19633\,(-), \\
a_{1,1}\lambda &= -0{,}00151\,58492\,(-), & l_1\lambda &= 0{,}01021\,16577\,(1), \\
\varepsilon &= \phantom{-}0{,}01021\,16577\,(-), & l_2\lambda^2 &= 0{,}00003\,00097\,(+), \\
a_{2,-2}\lambda^2 &= -0{,}00000\,01637\,(-), & l_3\lambda^3 &= 0{,}00000\,77466\,(-), \\
a_{2,0}\lambda^2 &= -0{,}00002\,39667\,(+), & \varepsilon' &= 0{,}00975\,60241\,(-), \\
a_{2,2}\lambda^2 &= -0{,}00000\,58793\,(+), & \delta &= 0{,}00985\,66771\,(-), \\
\varepsilon_1 &= \phantom{-}0{,}00003\,00097\,(+), & g &= 1{,}23269\,98742\,(-).
\end{aligned}$$

---

[8] Знак $(-)$ указывает здесь на то, что точное число по абсолютной величине менее написанного, знак $(+)$, что оно более написанного.

Обращаясь теперь уравнению (49) и ограничиваясь при вычислении наименьшего его корня восемью десятичными знаками, найдем

$$z = 0{,}01025064 \, (-).$$

Отсюда с таким же приближением выводящим

$$z - l_1\lambda = 0{,}00003\,898, \tag{51}$$

$$z - l_1\lambda - l_2\lambda^2 = 0{,}00000897, \tag{52}$$

$$z - l_1\lambda - l_2\lambda^2 - l_3\lambda^3 = 0{,}00000\,12^{3(-)}_{2(+)}. \tag{53}$$

Числа эти представляют высшие пределы погрешности, соответствующие случаям $n = 1$, $n = 2$ и $n = 3$.

Сопоставляя число (51) с величиной $\varepsilon_1$, которой может достигать совокупность членов второго порядка в выражении функции $\xi$ [формула (30)], приходим к заключено, что число это можно считать довольно удовлетворительным по точности высшим пределом погрешности.

Чтобы судить о степени точности высшего предела (52), необходимо вычислить дли выражений (30) по крайней мере коэффициенты четьего порядка относительно $\lambda$, и составить себе некоторое представление о величине

$$\varepsilon_2 = \{|a_{3,-3} - a_{1,-1}| + |a_{1,1}| + |a_{3,3}|\}\lambda^3.$$

Обращаясь к формулам (31) и (33), находим

$$a_{3,\pm 3} = -\alpha_{\pm 3,0}, \quad a_{3,\pm 1} = -a_{2,0}\alpha_{1,\pm 1} - \alpha_{\pm 1,1},$$

Мы должны, таким образом, вычислить

$$a_{2,0}a_{1,\pm 1}\lambda^3, \quad \alpha_{\pm 1,1}\lambda^3, \quad \text{и} \quad \alpha_{\pm 3,0}\lambda^3.$$

Из полученных выше численных результатов выводим

$$\begin{aligned} a_{2,0}a_{1,-1}\lambda^3 &= -0{,}00000\,02085\,(-), \\ a_{2,0}a_{1,1}\lambda^3 &= \phantom{-}0{,}00000\,00414\,(-). \end{aligned}$$

Для получения же $\alpha_{1,\pm 1}\lambda^3$ и $\alpha_{\pm 3,0}\lambda^3$ обращаемся к численным выводам Хилла, из которых заимствуем

$$\begin{aligned} \alpha_{-1,1}\lambda^3 &= \phantom{-}0{,}00000\,00616\,(-), \\ \alpha_{1,1}\lambda^3 &= -0{,}00000\,01417\,(-), \\ \alpha_{-3,0}\lambda^3 &= \phantom{-}0{,}00000\,00025\,(-), \\ \alpha_{3,0}\lambda^3 &= \phantom{-}0{,}00000\,00300\,(+) \end{aligned}$$

Отсюда находим

$$\begin{aligned} a_{3,-1}\lambda^3 &= \phantom{-}0{,}00000\,01469, \\ a_{3,1}\lambda^3 &= \phantom{-}0{,}00000\,01003, \\ a_{3,-3}\lambda^3 &= -0{,}00000\,00025\,(-), \\ a_{3,3}\lambda^3 &= -0{,}00000\,00300\,(+) \end{aligned}$$

и, наконец,

$$\varepsilon_2 = 0{,}00000\ 0280\ (-).$$

Таким образом

$$\varepsilon_2 < 0{,}00000028,$$

а число (53) не превосходит

$$0{,}00000\ 123.$$

Поэтому совокупности членов выше второго порядка в рядах (30) наверно будут численно менее

$$0{,}00000\ 151, \qquad (54)$$

и, следовательно, наш высший предел (52), почти в 6 раз больший этого числа, нельзя назвать удовлетворительным.

Еще менее удовлетворительным окажется, конечно, высший предел (53).

В заключение воспользуемся числами для опенки степени точности некоторых приближенных выражений $\xi$ и $\eta$.

Допустим, что в рядах (30) пренебрегается всеми членами выше первого порядка относительно $\lambda$ и что коэффициенты членов первого порядка вычисляются с пятью десятичными знаками, что дает

$$\xi = -0{,}00718 \cos 2\tau, \quad \eta = 0{,}01021 \sin 2\tau.$$

Принимая за высший предел погрешности для случая $n = 2$ число (54) и складывая его с $\varepsilon_1$, получим для высшего предела погрешности в случае $n = 1$ следующее число:

$$0{,}00003152.$$

Поэтому, замечая, что погрешности в коэффициентах наших приближенных выражений $\xi, \eta$, не достигают $0{,}000002$, заключаем, что вычисление этих выражений, каково бы ни было вещественное число $\tau$, может доставить для $\xi$ и $\eta$ величины, погрешности которых не будут превосходить $0{,}00004$.

Желая достигнуть большей точности при коэффициентах, вычисляемых с пятью десятичными знаками, мы должны принять в расчет члены второго порядка. При этом, останавливаясь на формулах

$$\begin{aligned}\xi &= 0{,}00002 & \eta &= 0{,}01021 \sin 2\tau \\ &\quad -0{,}00718 \cos 2\tau & &\quad +0{,}00001 \sin 4\tau, \\ &\quad +0{,}00001 \cos 4\tau \end{aligned}$$

можем рассчитывать получить для $\xi$ и $\eta$, величины, погрешности которых не будут превышать $0{,}00001$.

————